\documentclass[11pt,reqno]{amsart}
\usepackage{amsmath,amssymb}
\usepackage[colorlinks=true,linkcolor=magenta,citecolor=blue]{hyperref}
\usepackage[margin=1in]{geometry}
\usepackage{enumitem} 
\usepackage{mathabx}
\usepackage{epsfig}

\newtheorem{dn}{Definition}[section]
\newtheorem{dl}{Theorem}[section]
\newtheorem{md}{Proposition}[section]
\newtheorem{bd}{Lemma}[section]
\newtheorem{hq}{Corollary}[section]
\newtheorem{nx}{Remark}[section]
\newtheorem{vd}{Example}[section]
\newcommand{\R}{\mathbb{R}}

\newcommand{\ity}{\infty}

\newcommand{\bbd}{\begin{bd}}
\newcommand{\ebd}{\end{bd}}
\newcommand{\bdn}{\begin{dn}}
\newcommand{\edn}{\end{dn}}
\newcommand{\bhq}{\begin{hq}}
\newcommand{\ehq}{\end{hq}}
\newcommand{\bdl}{\begin{dl}}
\newcommand{\edl}{\end{dl}}
\newcommand{\bnx}{\begin{nx}}
\newcommand{\enx}{\end{nx}}
\newcommand{\bmd}{\begin{md}}
\newcommand{\emd}{\end{md}}
\newcommand{\bvd}{\begin{vd}}
\newcommand{\evd}{\end{vd}}

\title[Asymptotic profile of solutions to structurally damped $\sigma$-evolution equations]{Some remarks on the asymptotic profile of solutions to structurally damped $\sigma$-evolution equations}

\author{Tuan Anh Dao}
\address{Tuan Anh Dao \hfill\break
School of Applied Mathematics and Informatics, Hanoi University of Science and Technology, No.1 Dai Co Viet road, Hanoi, Vietnam \hfill\break
Faculty for Mathematics and Computer Science, TU Bergakademie Freiberg, Pr\"{u}ferstr. 9, 09596, Freiberg, Germany}
\email{anh.daotuan@hust.edu.vn}

\begin{document}
	
\begin{abstract}
In this paper, we are interested in analyzing the asymptotic profiles of solutions to the Cauchy problem for linear structurally damped $\sigma$-evolution equations in $L^2$-sense. Depending on the parameters $\sigma$ and $\delta$ we would like to not only indicate approximation formula of solutions but also recognize the optimality of their decay rates as well in the distinct cases of parabolic like damping and $\sigma$-evolution like damping. Moreover, such results are also discussed when we mix these two kinds of damping terms in a $\sigma$-evolution equation to investigate how each of them affects the asymptotic profile of solutions.
\end{abstract}
\maketitle

\section{Introduction}
Let us consider the following Cauchy problem for structurally damped $\sigma$-evolution equations:
\begin{equation}
\begin{cases}
u_{tt}+ (-\Delta)^\sigma u+ a(-\Delta)^{\delta_1} u_t+ b(-\Delta)^{\delta_2} u_t= 0, \\
u(0,x)= u_0(x),\quad u_t(0,x)=u_1(x),
\end{cases}
\label{equation1.1}
\end{equation}
where $\sigma \ge 1$, $0< \delta_1< \sigma/2< \delta_2< \sigma$ and $a,b=0,1$ with $(a,b) \neq (0,0)$.

At first, we recall some recent results related to the study of (\ref{equation1.1}) with $\sigma=1$ and $(a,b)=(1,0)$ or $(a,b)=(0,1)$, the so-called structurally damped wave equations, in the following form:
\begin{equation}
\begin{cases}
u_{tt}- \Delta u+ \mu (-\Delta)^{\delta} u_t= 0, \\
u(0,x)= u_0(x),\quad u_t(0,x)=u_1(x),
\end{cases}
\label{equation1.2}
\end{equation}
with $\delta \in (0,1)$ and $\mu>0$. The authors in \cite{DabbiccoReissig} succeeded in obtaining some of sharp $(L^1\cap L^2)- L^2$ estimates to (\ref{equation1.2}), i.e. the mixing of additional $L^1$ regularity for the data on the basis of $L^2- L^2$ estimates. A direct application of these estimates is to prove the global (in time) existence of small data energy solutions in low space dimensions to the corresponding semi-linear structurally damped wave models with power nonlinearties $|u|^p$. The point in the cited paper is that from the view of decay estimates they proposed to separate (\ref{equation1.2}) into ``parabolic like models" with $\delta \in (0,\frac{1}{2})$ and ``hyperbolic like models" with $\delta \in (\frac{1}{2},1)$. This comes from the properties of solutions to (\ref{equation1.2}) which change completely from the former case to the latter case. After that, in the quite recent paper \cite{IkehataTakeda} the asymptotic profile of solutions to (\ref{equation1.2}) and some of their optimal decay estimates have been explored so well. In particular, the authors provided the different approximation formulas of solutions by a constant multiple of a special function for large $t\ge 1$ corresponding to the cases $\delta \in (0,\frac{1}{2})$, $\delta= \frac{1}{2}$ and $\delta \in (\frac{1}{2},1)$.

Concerning the more general cases of $\sigma \ge 1$ to (\ref{equation1.1}), the following so-called structurally damped $\sigma$-evolution equations have been well-studied in several recent papers (see, for example, \cite{Duong2019,DabbiccoEbert2016,DabbiccoEbert2017,DuongKainaneReissig,DaoReissig1,DaoReissig2}):
\begin{equation}
\begin{cases}
u_{tt}+ (-\Delta)^\sigma u+ \mu (-\Delta)^\delta u_t= 0, \\
u(0,x)= u_0(x),\quad u_t(0,x)=u_1(x),
\end{cases}
\label{equation1.3}
\end{equation}
where $\sigma \ge 1$, $\delta \in (0,\sigma)$ and $\mu>0$. Namely, the results for decay rates of solutions to (\ref{equation1.3}) in the $L^2-L^2$ theory by assuming additional $L^1$ regularity for the data were derived in \cite{DuongKainaneReissig}. Quite recently, taking into considerations some of decay estimates for solutions to (\ref{equation1.3}) basing on the $L^q-L^q$ theory for any $q\in (1,\ity)$ the authors in \cite{DabbiccoEbert2017,DaoReissig1,DaoReissig2} have investigated $(L^m \cap L^q)- L^q$ and $L^q- L^q$ estimates with $q\in (1,\ity)$ and $m\in [1,q)$. More in detail, to establish this, they applied two main strategies including the theory of modified Bessel functions combined with Fa\`{a} di Bruno's formula and the Mikhlin-H\"{o}rmander multiplier theorem. By using the obtained decay estimates, the novelty of the cited papers are to prove the global (in time) existence of small data Sobolev solutions from suitable function spaces basing on $L^q$ spaces and to determine critical exponents as well to some semi-linear models with power nonlinearties $|u|^p$ or $|u_t|^p$. However, one may realizes that the asymptotic profiles of solutions and the optimality of their decay rates have not been indicated in the above mentioned references clearly. For this reason, one of the main goals of this paper is to report such results for solutions to (\ref{equation1.3}).

According the classification of (\ref{equation1.3}) proposed in \cite{DabbiccoEbert2016,DabbiccoEbert2017,DaoReissig1,DaoReissig2}, here we want to distinguish (\ref{equation1.1}) into three main models depending on the parameters $a$ and $b$. In particular, the first model of our considerations is the $\sigma$-evolution equations with parabolic like structural damping corresponding to the case $(a,b)=(1,0)$. In this model, we are going to show that the asymptotic profile of solutions to (\ref{equation1.1}) is the same as that to the following anomalous diffusion equations (see later, Theorem \ref{theorem1.1}):
\begin{equation}
v_t+ (-\Delta)^{\sigma-\delta_1} v= 0, \qquad v(0,x)= v_0(x),
\label{equation1.4}
\end{equation}
for a suitable choice of data $v_0$. The second one is the model with $\sigma$-evolution like structural damping corresponding to the case $(a,b)=(0,1)$, the so-called ``hyperbolic like models" in the case $\sigma=1$. We recognize that some kind of wave structure appears and oscillations come into play from the asymptotic profile of solutions in this model (see later, Theorem \ref{theorem1.2}). This means the above mentioned diffusion phenomenon does not happen. Our interest is the last model with mixing two distinct kinds of structural damping including parabolic like damping and $\sigma$-evolution like damping corresponding to the case $(a,b)=(1,1)$, the so-called double damping terms (see, for instance, \cite{DaoMichihisa,IkehataMichihisa,IkehataSawada}). This connection brings some interesting properties for solutions to (\ref{equation1.1}) in the case $\delta_1+ \delta_2> \sigma$ which inherit from the two former models (see later, Theorem \ref{theorem1.3}). More precisely, by the presence of parabolic like damping, on the one hand, the asymptotic profile of solutions to (\ref{equation1.1}) is also the same as that to (\ref{equation1.4}). On the other hand, the solutions to (\ref{equation1.1}) possess the same regularity as that to the second model by the presence of $\sigma$-evolution like damping. Analyzing these properties is the second main goal of the this paper.

\subsection{Notations}
\begin{itemize}[leftmargin=*]
\item We write $f\lesssim g$ when there exists a constant $C>0$ such that $f\le Cg$, and $f \approx g$ when $g\lesssim f\lesssim g$.
\item As usual, $H^a$ and $\dot{H}^a$, with $a \ge 0$, denote Bessel and Riesz potential spaces based on $L^2$ spaces. Here $\big<D\big>^a$ and $|D|^a$ stand for the pseudo-differential operators with symbols $\big<\xi\big>^a$ and $|\xi|^a$, respectively.
\item We denote $\widehat{w}(t,\xi):= \mathfrak{F}_{x\rightarrow \xi}\big(w(t,x)\big)$ as the Fourier transform with respect to the space variable of a function $w= w(t,x)$. 
\item We put $[s]^+:= \max\{s,0\}$ as the positive part of $s \in \R$. Moreover, we fix the constant $m_0:=\frac{2m}{2-m}$, that is, $\frac{1}{m_0}=\frac{1}{m}- \frac{1}{2}$ with $m \in [1,2)$.
\item Let $\chi= \chi(|\xi|)$ be a $\mathcal{C}_0^\ity(\R^n)$ cut-off nonnegative function equal to $1$ for small $|\xi|$ and vanishing for large $|\xi|$. We decompose a function $w= w(t,x)$ into two parts localized separately to low and high frequencies as follows:
$$ w(t,x)= w_{\text{\fontshape{n}\selectfont low}}(t,x)+ w_{\text{\fontshape{n}\selectfont high}}(t,x), $$
where
$$w_{\text{\fontshape{n}\selectfont low}}(t,x)= \mathfrak{F}^{-1}\big(\chi(|\xi|)\widehat{w}(t,\xi)\big)\quad \text{ and }\quad w_{\text{\fontshape{n}\selectfont high}}(t,x)= \mathfrak{F}^{-1}\Big(\big(1-\chi(|\xi|)\big)\widehat{w}(t,\xi)\Big). $$
\item For later convenience, we denote the following quantity:
$$P_1:= \int_{\R^n}u_1(x)dx. $$
\end{itemize}

\subsection{Main results}
The following results describe the large time behavior of solutions to (\ref{equation1.1}).
\bdl[\textbf{$a=1$ and $b=0$}] \label{theorem1.1}
Let $j=0,1$ and $s \ge 0$. We assume the condition $n>4\delta_1$ and the initial data
$$ (u_0,u_1) \in \mathcal{A}^1_0 \times \mathcal{A}^1_1:= (L^1 \cap H^{s+j\sigma}) \times (L^1 \cap H^{[s+(j-1)\sigma]^+}). $$
Then, the Sobolev solutions to (\ref{equation1.1}) satisfy the following estimates for large $t \ge 1$:
\begin{equation}
\Big\|\partial_t^j |D|^s \Big(u(t,\cdot)- P_1\,\mathfrak{F}^{-1}\Big(\dfrac{e^{-t|\xi|^{2(\sigma-\delta_1)}}}{|\xi|^{2\delta_1}}\Big)(t,\cdot)\Big)\Big\|_{L^2}= o\Big(t^{-\frac{n}{4(\sigma-\delta_1)}-\frac{s}{2(\sigma-\delta_1)}-j+ \frac{\delta_1}{\sigma-\delta_1}}\Big). \label{theorem1.1.1}
\end{equation}
Moreover, if $P_1 \neq 0$, then the following estimates hold for large $t \ge 1$:
\begin{equation}
C_1\,t^{-\frac{n}{4(\sigma-\delta_1)}-\frac{s}{2(\sigma-\delta_1)}-j+ \frac{\delta_1}{\sigma-\delta_1}}\le \big\|\partial_t^j |D|^s  u(t,\cdot)\big\|_{L^2}\le C_2\,t^{-\frac{n}{4(\sigma-\delta_1)}-\frac{s}{2(\sigma-\delta_1)}-j+ \frac{\delta_1}{\sigma-\delta_1}}, \label{theorem1.1.2}
\end{equation}
where $C_1$ and $C_2$ are some suitable positive constants.
\edl

\bdl[\textbf{$a=0$ and $b=1$}] \label{theorem1.2}
Let $j=0,1$ and $s \ge 0$. We assume the condition $n>2\sigma$ and the initial data
$$ (u_0,u_1) \in \mathcal{A}^2_0 \times \mathcal{A}^2_1:= (L^1 \cap H^{s+2j\delta_2}) \times (L^1 \cap H^{[s+2(j-1)\delta_2]^+}). $$
Then, the Sobolev solutions to (\ref{equation1.1}) satisfy the following estimates for large $t \ge 1$:
\begin{align}
\Big\||D|^s \Big(u(t,\cdot)- P_1\,\mathfrak{F}^{-1}\Big(e^{-\frac{1}{2}t|\xi|^{2\delta_2}}\dfrac{\sin(t|\xi|^\sigma)}{|\xi|^\sigma}\Big)(t,\cdot)\Big)\Big\|_{L^2}&= o\Big(t^{-\frac{n}{4\delta_2}-\frac{s}{2\delta_2}+\frac{\sigma}{2\delta_2}}\Big), \label{theorem1.2.1} \\
\Big\|\partial_t |D|^s u(t,\cdot)- P_1\,|D|^s \mathfrak{F}^{-1}\Big(e^{-\frac{1}{2}t|\xi|^{2\delta_2}}\cos(t|\xi|^\sigma)\Big)(t,\cdot)\Big\|_{L^2}&= o\Big(t^{-\frac{n}{4\delta_2}-\frac{s}{2\delta_2}}\Big). \label{theorem1.2.2}
\end{align}
Moreover, if $P_1 \neq 0$, then the following estimates hold for large $t \ge 1$:
\begin{equation}
C_1\,t^{-\frac{n}{4\delta_2}-\frac{s+(j-1)\sigma}{2\delta_2}}\le \big\|\partial_t^j |D|^s  u(t,\cdot)\big\|_{L^2}\le C_2\,t^{-\frac{n}{4\delta_2}-\frac{s+(j-1)\sigma}{2\delta_2}}, \label{theorem1.2.3}
\end{equation}
where $C_1$ and $C_2$ are some suitable positive constants.
\edl

\bdl[\textbf{$a=1$ and $b=1$}] \label{theorem1.3}
Let $j=0,1$ and $s \ge 0$. Let us assume $\delta_1+\delta_2>\sigma$. We suppose the condition $n>4\delta_1$ and the initial data
$$ (u_0,u_1) \in \mathcal{A}^3_0 \times \mathcal{A}^3_1:= (L^1 \cap H^{s+2j\delta_2}) \times (L^1 \cap H^{[s+2(j-1)\delta_2]^+}). $$
Then, the Sobolev solutions to (\ref{equation1.1}) satisfy the following estimates for large $t \ge 1$:
\begin{equation}
\Big\|\partial_t^j |D|^s \Big(u(t,\cdot)- P_1\,\mathfrak{F}^{-1}\Big(\dfrac{e^{-t|\xi|^{2(\sigma-\delta_1)}}}{|\xi|^{2\delta_1}}\Big)(t,\cdot)\Big)\Big\|_{L^2}= o\Big(t^{-\frac{n}{4(\sigma-\delta_1)}-\frac{s}{2(\sigma-\delta_1)}-j+ \frac{\delta_1}{\sigma-\delta_1}}\Big). \label{theorem1.3.1}
\end{equation}
Moreover, if $P_1 \neq 0$, then the following estimates hold for large $t \ge 1$:
\begin{equation}
C_1\,t^{-\frac{n}{4(\sigma-\delta_1)}-\frac{s}{2(\sigma-\delta_1)}-j+ \frac{\delta_1}{\sigma-\delta_1}}\le \big\|\partial_t^j |D|^s  u(t,\cdot)\big\|_{L^2}\le C_2\,t^{-\frac{n}{4(\sigma-\delta_1)}-\frac{s}{2(\sigma-\delta_1)}-j+ \frac{\delta_1}{\sigma-\delta_1}}, \label{theorem1.3.2}
\end{equation}
where $C_1$ and $C_2$ are some suitable positive constants.
\edl

\textbf{The organization of this paper} is as follows: In Section \ref{Sec.2}, we present preliminary knowledge as the representation of solutions, pointwise estimates in Fourier space and some of decay estimates for solutions to (\ref{equation1.1}) in Sections \ref{Sec.2.1}, \ref{Sec.2.2} and \ref{Sec.2.3}, respectively. Then, we prove main results for the asymptotic profile of solutions to (\ref{equation1.1}) and indicate their optimal decay estimates as well in Section \ref{Sec.3}. In particular, the proofs in Sections \ref{Sec.3.1}, \ref{Sec.3.2} and \ref{Sec.3.3} correspond to the cases $(a,b)= (1,0)$, $(a,b)= (0,1)$ and $(a,b)= (1,1)$.

\section{Preliminaries} \label{Sec.2}

\subsection{Representation of solutions} \label{Sec.2.1}
At first, using partial Fourier transformation to (\ref{equation1.1}) we obtain the following Cauchy problem:
\begin{equation}
\widehat{u}_{tt}+ \big(a|\xi|^{2\delta_1}+ b|\xi|^{2\delta_2}\big) \widehat{u}_t+ |\xi|^{2\sigma} \widehat{u}=0,\quad \widehat{u}(0,\xi)= \widehat{u_0}(\xi),\quad \widehat{u}_t(0,\xi)= \widehat{u_1}(\xi). \label{equation2.1}
\end{equation}
The characteristic roots are
$$ \lambda_{1,2}=\lambda_{1,2}(\xi)= \frac{1}{2}\Big(-\big(a|\xi|^{2\delta_1}+ b|\xi|^{2\delta_2}\big)\pm \sqrt{\big(a|\xi|^{2\delta_1}+ b|\xi|^{2\delta_2}\big)^2-4|\xi|^{2\sigma}}\Big). $$
The solutions to (\ref{equation2.1}) are presented by the following formula (here we assume $\lambda_{1}\neq \lambda_{2}$):
\begin{equation*}
\widehat{u}(t,\xi)= \frac{\lambda_1 e^{\lambda_2 t}-\lambda_2 e^{\lambda_1 t}}{\lambda_1- \lambda_2}\widehat{u_0}(\xi)+ \frac{e^{\lambda_1 t}-e^{\lambda_2 t}}{\lambda_1- \lambda_2}\widehat{u_1}(\xi)=: \widehat{K_0}(t,\xi)\widehat{u_0}(\xi)+\widehat{K_1}(t,\xi)\widehat{u_1}(\xi).
\end{equation*}
Now depending on the parameters $a$ and $b$ we shall decompose the above representation formula of solutions to (\ref{equation2.1}) into several parts as follows:
\begin{itemize}[leftmargin=*]
\item $a=1$ and $b=0,\,1$:
$$ \widehat{u}(t,\xi)= \big(\widehat{K^1_0}(t,\xi)+ \widehat{K^2_0}(t,\xi)\big)\widehat{u_0}(\xi)+ \big(\widehat{K^1_1}(t,\xi)+ \widehat{K^2_1}(t,\xi)\big)\widehat{u_1}(\xi), $$
where
\begin{align*}
\widehat{K^1_0}(t,\xi)&= \frac{-\lambda_2 e^{\lambda_1 t}}{\lambda_1- \lambda_2}, \quad &\widehat{K^2_0}(t,\xi)&= \frac{\lambda_1 e^{\lambda_2 t}}{\lambda_1- \lambda_2}, \\ 
\widehat{K^1_1}(t,\xi)&= \frac{e^{\lambda_1 t}}{\lambda_1- \lambda_2}, \quad &\widehat{K^2_1}(t,\xi)&= \frac{-e^{\lambda_2 t}}{\lambda_1- \lambda_2}.
\end{align*}
\item $a=0$ and $b=1$:
$$ \widehat{u}(t,\xi)= \big(\widehat{K^{\cos}_0}(t,\xi)+ \widehat{K^{\sin}_0}(t,\xi)\big)\widehat{u_0}(\xi)+ \widehat{K_1}(t,\xi)\widehat{u_1}(\xi), $$
where
$$ \widehat{K^{\cos}_0}(t,\xi)= e^{-\frac{1}{2}t|\xi|^{2\delta_2}}\cos\big(t|\xi|^\sigma f(|\xi|)\big), \quad \widehat{K^{\sin}_0}(t,\xi)= e^{-\frac{1}{2}t|\xi|^{2\delta_2}}|\xi|^{2\delta_2} \dfrac{\sin\big(t|\xi|^\sigma f(|\xi|)\big)}{2|\xi|^\sigma f(|\xi|)}, $$
$$ \widehat{K_1}(t,\xi)= e^{-\frac{1}{2}t|\xi|^{2\delta_2}}\dfrac{\sin\big(t|\xi|^\sigma f(|\xi|)\big)}{|\xi|^\sigma f(|\xi|)}, $$
with
$$ f(|\xi|)=
\begin{cases}
\sqrt{1-\frac{1}{4}|\xi|^{4\delta_2- 2\sigma}}&\text{ for small }|\xi|, \\
i\sqrt{\frac{1}{4}|\xi|^{4\delta_2- 2\sigma}-1}&\text{ for large }|\xi|.
\end{cases} $$
\end{itemize}

\subsection{Pointwise estimates in Fourier space} \label{Sec.2.2}
Taking account of the cases of small and large frequencies separately we have the asymptotic behavior of the characteristic roots as follows:
\begin{align*}
&1.\,\, a=1 \text{ and }b=0:& &\lambda_1\sim -|\xi|^{2(\sigma- \delta_1)},\qquad \lambda_2\sim -|\xi|^{2\delta_1},\qquad \lambda_1-\lambda_2 \sim |\xi|^{2\delta_1} \quad \text{ for small } |\xi|, \\
& & &\text{and }\lambda_{1,2} \sim -|\xi|^{2\delta_1}\pm i|\xi|^\sigma,\qquad \lambda_1-\lambda_2 \sim i|\xi|^\sigma \quad \text{ for large } |\xi|, \\
&2.\,\, a=0 \text{ and }b=1:& &\lambda_{1,2} \sim -|\xi|^{2\delta_2}\pm i|\xi|^\sigma,\qquad \lambda_1-\lambda_2 \sim i|\xi|^\sigma \quad \text{ for small } |\xi|, \\
& & &\text{and }\lambda_1\sim -|\xi|^{2(\sigma- \delta_2)},\qquad \lambda_2\sim -|\xi|^{2\delta_2},\qquad \lambda_1-\lambda_2 \sim |\xi|^{2\delta_2} \quad \text{ for large } |\xi|, \\
&3.\,\, a=1 \text{ and }b=1:& &\lambda_1\sim -|\xi|^{2(\sigma- \delta_1)},\qquad \lambda_2\sim -|\xi|^{2\delta_1},\qquad \lambda_1-\lambda_2 \sim |\xi|^{2\delta_1} \quad \text{ for small } |\xi|, \\
& & &\text{and }\lambda_1\sim -|\xi|^{2(\sigma- \delta_2)},\qquad \lambda_2\sim -|\xi|^{2\delta_2},\qquad \lambda_1-\lambda_2 \sim |\xi|^{2\delta_2} \quad \text{ for large } |\xi|.
\end{align*}

\subsubsection{The case $a=1$ and $b=0$}
\bbd \label{lemma2.1}
Let $s\ge 0$ and $j=0,\,1$. Then, the following estimates hold:
\begin{align*}
|\xi|^s \chi(|\xi|)\big|\partial^j_t \widehat{K^1_0}(t,\xi)\big| &\lesssim e^{-c_0 t|\xi|^{2(\sigma-\delta_1)}}|\xi|^{s+2j(\sigma-\delta_1)}, \\
|\xi|^s \chi(|\xi|)\big|\partial^j_t \widehat{K^2_0}(t,\xi)\big| &\lesssim e^{-c_0 t|\xi|^{2\delta_1}}|\xi|^{s+2j\delta_1+2(\sigma-2\delta_1)}, \\
|\xi|^s \chi(|\xi|)\big|\partial^j_t \widehat{K^1_1}(t,\xi)\big| &\lesssim e^{-c_0 t|\xi|^{2(\sigma-\delta_1)}}|\xi|^{s+2j(\sigma-\delta_1)-2\delta_1}, \\
|\xi|^s \chi(|\xi|)\big|\partial^j_t \widehat{K^2_1}(t,\xi)\big| &\lesssim e^{-c_0 t|\xi|^{2\delta_1}}|\xi|^{s+2(j-1)\delta_1},
\end{align*}
and
\begin{align*}
|\xi|^s \big(1-\chi(|\xi|)\big)\big|\partial^j_t \widehat{K_0}(t,\xi)\big| &\lesssim e^{-c_0 t|\xi|^{2\delta_1}}|\xi|^{s+j\sigma}, \\
|\xi|^s \big(1-\chi(|\xi|)\big)\big|\partial^j_t \widehat{K_1}(t,\xi)\big| &\lesssim e^{-c_0 t|\xi|^{2\delta_1}}|\xi|^{s+(j-1)\sigma},
\end{align*}
where $c_0$ is a suitable positive constant.
\ebd

\subsubsection{The case $a=0$ and $b=1$}
\bbd \label{lemma2.2}
Let $s\ge 0$ and $j=0,\,1$. Then, the following estimates hold:
\begin{align*}
|\xi|^s \chi(|\xi|)\big|\partial^j_t \widehat{K^{\cos}_0}(t,\xi)\big| &\lesssim e^{-c_0 t|\xi|^{2\delta_2}}|\xi|^{s+j\sigma}, \\
|\xi|^s \chi(|\xi|)\big|\partial^j_t \widehat{K^{\sin}_0}(t,\xi)\big| &\lesssim e^{-c_0 t|\xi|^{2\delta_2}}|\xi|^{s+(j-1)\sigma+2\delta_2}, \\
|\xi|^s \chi(|\xi|)\big|\partial^j_t \widehat{K_1}(t,\xi)\big| &\lesssim e^{-c_0 t|\xi|^{2\delta_2}}|\xi|^{s+(j-1)\sigma},
\end{align*}
and
\begin{align*}
|\xi|^s \big(1-\chi(|\xi|)\big)\big|\partial^j_t \widehat{K_0}(t,\xi)\big| &\lesssim e^{-c_0 t|\xi|^{2(\sigma-\delta_2)}}|\xi|^{s+2j(\sigma-\delta_2)}+ e^{-c_0 t|\xi|^{2\delta_2}}|\xi|^{s+2j\delta+2(\sigma-2\delta_2)}, \\
|\xi|^s \big(1-\chi(|\xi|)\big)\big|\partial^j_t \widehat{K_1}(t,\xi)\big| &\lesssim e^{-c_0 t|\xi|^{2(\sigma-\delta_2)}}|\xi|^{s+2j(\sigma-\delta_2)-2\delta_2}+ e^{-c_0 t|\xi|^{2\delta_2}}|\xi|^{s+2(j-1)\delta_2},
\end{align*}
where $c_0$ is a suitable positive constant.
\ebd

\subsubsection{The case $a=1$ and $b=1$}
\bbd \label{lemma2.3}
Let $s\ge 0$ and $j=0,\,1$. Then, the following estimates hold:
\begin{align*}
|\xi|^s \chi(|\xi|)\big|\partial^j_t \widehat{K^1_0}(t,\xi)\big| &\lesssim e^{-c_0 t|\xi|^{2(\sigma-\delta_1)}}|\xi|^{s+2j(\sigma-\delta_1)}, \\
|\xi|^s \chi(|\xi|)\big|\partial^j_t \widehat{K^2_0}(t,\xi)\big| &\lesssim e^{-c_0 t|\xi|^{2\delta_1}}|\xi|^{s+2j\delta_1+2(\sigma-2\delta_1)}, \\
|\xi|^s \chi(|\xi|)\big|\partial^j_t \widehat{K^1_1}(t,\xi)\big| &\lesssim e^{-c_0 t|\xi|^{2(\sigma-\delta_1)}}|\xi|^{s+2j(\sigma-\delta_1)-2\delta_1}, \\
|\xi|^s \chi(|\xi|)\big|\partial^j_t \widehat{K^2_1}(t,\xi)\big| &\lesssim e^{-c_0 t|\xi|^{2\delta_1}}|\xi|^{s+2(j-1)\delta_1},
\end{align*}
and
\begin{align*}
|\xi|^s \big(1-\chi(|\xi|)\big)\big|\partial^j_t \widehat{K_0}(t,\xi)\big| &\lesssim e^{-c_0 t|\xi|^{2(\sigma-\delta_2)}}|\xi|^{s+2j(\sigma-\delta_2)}+ e^{-c_0 t|\xi|^{2\delta_2}}|\xi|^{s+2j\delta_2+2(\sigma-2\delta_2)}, \\
|\xi|^s \big(1-\chi(|\xi|)\big)\big|\partial^j_t \widehat{K_1}(t,\xi)\big| &\lesssim e^{-c_0 t|\xi|^{2(\sigma-\delta_2)}}|\xi|^{s+2j(\sigma-\delta_2)-2\delta_2}+ e^{-c_0 t|\xi|^{2\delta_2}}|\xi|^{s+2(j-1)\delta_2},
\end{align*}
where $c_0$ is a suitable positive constant.
\ebd

\subsection{Decay estimates} \label{Sec.2.3}
\subsubsection{The case $a=1$ and $b=0$}
\bbd \label{lemma2.4}
Let $s\ge 0$ and $j=0,\,1$. Then, the following estimates hold for $m \in [1,2)$:
\begin{align}
\big\|\partial_t^j |D|^s \big(K^1_{0\text{\fontshape{n}\selectfont low}}(t,x) \ast u_0(x)\big)(t,\cdot)\big\|_{L^2}&\lesssim (1+t)^{-\frac{n}{2(\sigma-\delta_1)}(\frac{1}{m}-\frac{1}{2})-\frac{s}{2(\sigma-\delta_1)}-j}\|u_0\|_{L^m}, \label{lemma2.4.1} \\
\big\|\partial_t^j |D|^s \big(K^2_{0\text{\fontshape{n}\selectfont low}}(t,x) \ast u_0(x)\big)(t,\cdot)\big\|_{L^2}&\lesssim (1+t)^{-\frac{n}{2\delta_1}(\frac{1}{m}-\frac{1}{2})-\frac{s}{2\delta_1}-j-\frac{\sigma-2\delta_1}{\delta_1}}\|u_0\|_{L^m}, \label{lemma2.4.2}
\end{align}
for any space dimensions $n\ge 1$, and
\begin{align}
\big\|\partial_t^j |D|^s \big(K^1_{1\text{\fontshape{n}\selectfont low}}(t,x) \ast u_1(x)\big)(t,\cdot)\big\|_{L^2}&\lesssim (1+t)^{-\frac{n}{2(\sigma-\delta_1)}(\frac{1}{m}-\frac{1}{2})-\frac{s}{2(\sigma-\delta_1)}-j+\frac{\delta_1}{\sigma-\delta_1}}\|u_1\|_{L^m}, \label{lemma2.4.3} \\
\big\|\partial_t^j |D|^s \big(K^2_{1\text{\fontshape{n}\selectfont low}}(t,x) \ast u_1(x)\big)(t,\cdot)\big\|_{L^2}&\lesssim (1+t)^{-\frac{n}{2\delta_1}(\frac{1}{m}-\frac{1}{2})-\frac{s}{2\delta_1}-j+1}\|u_1\|_{L^m}, \label{lemma2.4.4}
\end{align}
for any space dimensions $n> 2m_0\delta_1$. Moreover, the following estimates hold for $m \in [1,2]$:
\begin{align}
\big\|\partial_t^j |D|^s \big(K_{0\text{\fontshape{n}\selectfont high}}(t,x) \ast u_0(x)\big)(t,\cdot)\big\|_{L^2}&\lesssim e^{-ct}t^{-\frac{n}{2\delta_1}(\frac{1}{m}-\frac{1}{2})-\frac{s+j\sigma-a}{2\delta_1}}\|u_0\|_{H^a_m}, \label{lemma2.4.5} \\
\big\|\partial_t^j |D|^s \big(K_{1\text{\fontshape{n}\selectfont high}}(t,x) \ast u_1(x)\big)(t,\cdot)\big\|_{L^2}&\lesssim e^{-ct}t^{-\frac{n}{2\delta_1}(\frac{1}{m}-\frac{1}{2})-\frac{s+j\sigma-a}{2\delta_1}}\|u_1\|_{H^{[a-\sigma]^+}_m}, \label{lemma2.4.6}
\end{align}
for all $a\ge 0$ and for any space dimensions $n\ge 1$, where $c$ is a suitable positive constant.
\ebd
\begin{proof}
First, we shall prove (\ref{lemma2.4.1}). By the first estimate in Lemma \ref{lemma2.1}, we apply Parseval-Plancherel formula and H\"{o}lder's inequality to obtain the following estimate:
\begin{align*}
\big\|\partial_t^j |D|^s \big(K^1_{0\text{\fontshape{n}\selectfont low}}(t,x) \ast u_0(x)\big)(t,\cdot)\big\|_{L^2}&= \big\| |\xi|^s \chi(|\xi|)\partial^j_t \widehat{K^1_0}(t,\xi)\widehat{u_0}(\xi)\big\|_{L^2} \\ 
&\lesssim \big\| e^{-ct|\xi|^{2(\sigma-\delta_1)}}|\xi|^{s+2j(\sigma-\delta_1)}\big\|_{L^{m_0}}\, \|\widehat{u_0}\|_{L^{m'}}.
\end{align*}
Thanks to the Hausdorff-Young inequality, we can control $\|\widehat{u_0}\|_{L^{m'}}$ by $\|u_0\|_{L^m}$. Hence, we have only to control the $L^{m_0}$ norm of the above multiplier. Using Lemma \ref{LemmaL1normEstimate} gives
$$ \big\| e^{-ct|\xi|^{2(\sigma-\delta_1)}}|\xi|^{s+2j(\sigma-\delta_1)}\big\|_{L^{m_0}}\lesssim (1+t)^{-\frac{n}{2m_0(\sigma-\delta_1)}-\frac{s}{2(\sigma-\delta_1)}-j}= (1+t)^{-\frac{n}{2(\sigma-\delta_1)}(\frac{1}{m}-\frac{1}{2})-\frac{s}{2(\sigma-\delta_1)}-j}. $$
This completes the proof of (\ref{lemma2.4.1}). In the same way we may arrive at the estimates from (\ref{lemma2.4.2}) to (\ref{lemma2.4.4}). Then, in order to indicate (\ref{lemma2.4.5}) and (\ref{lemma2.4.6}), we repeat the proof of (\ref{lemma2.4.1}) by using a suitable regularity of the data $u_0$ and $u_1$. Indeed, by Lemma \ref{lemma2.1} we get
\begin{align*}
\big\|\partial_t^j |D|^s \big(K_{0\text{\fontshape{n}\selectfont high}}(t,x) \ast u_0(x)\big)(t,\cdot)\big\|_{L^2}&= \big\| |\xi|^s \big(1-\chi(|\xi|)\big)\partial^j_t \widehat{K_0}(t,\xi)\widehat{u_0}(\xi)\big\|_{L^2} \\ 
&\lesssim \big\|e^{-c_0 t|\xi|^{2\delta_1}}|\xi|^{s+j\sigma}\widehat{u_0}(\xi)\big\|_{L^2} \lesssim e^{-\frac{c_0}{2}t}\big\|e^{-\frac{c_0}{2}t|\xi|^{2\delta_1}}|\xi|^{s+j\sigma}\widehat{u_0}(\xi)\big\|_{L^2} \\
&\lesssim e^{-ct}\big\|e^{-\frac{c_0}{2}t|\xi|^{2\delta_1}}|\xi|^{s+j\sigma-a}\big\|_{L^{m_0}}\big\| |\xi|^a\widehat{u_0}(\xi)\big\|_{L^{m'}},\,\, \text{ where }c:= \frac{c_0}{2} \\
&\lesssim e^{-ct}t^{-\frac{n}{2\delta_1}(\frac{1}{m}-\frac{1}{2})-\frac{s+j\sigma-a}{2\delta_1}}\|u_0\|_{H^a_m}.
\end{align*}
This completes the proof of (\ref{lemma2.4.5}). By an analogous argument we may also conclude (\ref{lemma2.4.6}). Therefore, Lemma \ref{lemma2.4} is proved.
\end{proof}

\noindent Hence, we obtain decay properties of Sobolev solutions to (\ref{equation1.1}).
\bmd \label{proposition2.1}
Let $s \ge 0$ and $j=0,1$. Let $m \in [1,2)$. We assume the condition $n>2m_0\delta_1$ and the initial data
$$ (u_0,u_1) \in (L^m \cap H^{s+j\sigma}) \times (L^m \cap H^{[s+(j-1)\sigma]^+}). $$
Then, the Sobolev solutions to (\ref{equation1.1}) satisfy the following $(L^m \cap L^2)-L^2$ estimates:
\begin{align*}
\big\|\partial_t^j |D|^s u(t,\cdot)\big\|_{L^2} &\lesssim (1+t)^{-\frac{n}{2(\sigma-\delta_1)}(\frac{1}{m}-\frac{1}{2})-\frac{s}{2(\sigma-\delta_1)}-j}\|u_0\|_{L^m \cap H^{s+j\sigma}} \\ 
&\qquad + (1+t)^{-\frac{n}{2(\sigma-\delta_1)}(\frac{1}{m}-\frac{1}{2})-\frac{s}{2(\sigma-\delta_1)}-j+ \frac{\delta_1}{\sigma-\delta_1}}\|u_1\|_{L^m \cap H^{[s+(j-1)\sigma]^+}}.
\end{align*}
\emd

\subsubsection{The case $a=0$ and $b=1$}
\bbd \label{lemma2.5}
Let $s\ge 0$ and $j=0,\,1$. Then, the following estimates hold for $m \in [1,2)$:
\begin{align}
\big\|\partial_t^j |D|^s \big(K^{\cos}_{0\text{\fontshape{n}\selectfont low}}(t,x) \ast u_0(x)\big)(t,\cdot)\big\|_{L^2}&\lesssim (1+t)^{-\frac{n}{2\delta_2}(\frac{1}{m}-\frac{1}{2})-\frac{s+j\sigma}{2\delta_2}}\|u_0\|_{L^m}, \label{lemma2.5.1} \\
\big\|\partial_t^j |D|^s \big(K^{\sin}_{0\text{\fontshape{n}\selectfont low}}(t,x) \ast u_0(x)\big)(t,\cdot)\big\|_{L^2}&\lesssim (1+t)^{-\frac{n}{2\delta_2}(\frac{1}{m}-\frac{1}{2})-\frac{s+(j-1)\sigma}{2\delta_2}-1}\|u_0\|_{L^m}, \label{lemma2.5.2}
\end{align}
for any space dimensions $n\ge 1$, and
\begin{equation}
\big\|\partial_t^j |D|^s \big(K_{1\text{\fontshape{n}\selectfont low}}(t,x) \ast u_1(x)\big)(t,\cdot)\big\|_{L^2}\lesssim (1+t)^{-\frac{n}{2\delta_2}(\frac{1}{m}-\frac{1}{2})-\frac{s+(j-1)\sigma}{2\delta_2}}\|u_1\|_{L^m}, \label{lemma2.5.3}
\end{equation}
for any space dimensions $n> m_0\sigma$. Moreover, the following estimates hold for $m \in [1,2]$:
\begin{align}
&\big\|\partial_t^j |D|^s \big(K_{0\text{\fontshape{n}\selectfont high}}(t,x) \ast u_0(x)\big)(t,\cdot)\big\|_{L^2} \nonumber \\
&\qquad \lesssim e^{-ct}\Big(t^{-\frac{n}{2(\sigma-\delta_2)}(\frac{1}{m}-\frac{1}{2})-\frac{s-a}{2(\sigma-\delta_2)}-j}+ t^{-\frac{n}{2\delta_2}(\frac{1}{m}-\frac{1}{2})-\frac{s-a}{2\delta_2}-j+ \frac{2\delta_2-\sigma}{\delta_2}}\Big)\|u_0\|_{H^a_m}, \label{lemma2.5.4}\\
&\big\|\partial_t^j |D|^s \big(K_{1\text{\fontshape{n}\selectfont high}}(t,x) \ast u_1(x)\big)(t,\cdot)\big\|_{L^2} \nonumber \\
&\qquad \lesssim e^{-ct}\Big(t^{-\frac{n}{2(\sigma-\delta_2)}(\frac{1}{m}-\frac{1}{2})-\frac{s-a}{2(\sigma-\delta_2)}-j}+ t^{-\frac{n}{2\delta_2}(\frac{1}{m}-\frac{1}{2})-\frac{s-a}{2\delta_2}-j}\Big)\|u_0\|_{H^{[a-2\delta_2]^+}_m}, \label{lemma2.5.5}
\end{align}
for all $a\ge 0$ and for any space dimensions $n\ge 1$, where $c$ is a suitable positive constant.
\ebd
\begin{proof}
Following the proof of Lemma \ref{lemma2.4} we may conclude all the desired estimates in Lemma \ref{lemma2.5} by the aid of the statements in Lemma \ref{lemma2.2}.
\end{proof}

\noindent Hence, we obtain decay properties of Sobolev solutions to (\ref{equation1.1}).
\bmd \label{proposition2.2}
Let $s \ge 0$ and $j=0,1$. Let $m \in [1,2)$. We assume the condition $n>m_0\sigma$ and the initial data
$$ (u_0,u_1) \in (L^m \cap H^{s+2j\delta_2}) \times (L^m \cap H^{[s+2(j-1)\delta_2]^+}). $$
Then, the Sobolev solutions to (\ref{equation1.1}) satisfy the following $(L^m \cap L^2)-L^2$ estimates:
\begin{align*}
\big\|\partial_t^j |D|^s u(t,\cdot)\big\|_{L^2} &\lesssim (1+t)^{-\frac{n}{2\delta_2}(\frac{1}{m}-\frac{1}{2})-\frac{s+j\sigma}{2\delta_2}}\|u_0\|_{L^m \cap H^{s+2j\delta_2}} \\ 
&\qquad + (1+t)^{-\frac{n}{2\delta_2}(\frac{1}{m}-\frac{1}{2})-\frac{s+(j-1)\sigma}{2\delta_2}}\|u_1\|_{L^m \cap H^{[s+2(j-1)\delta_2]^+}}.
\end{align*}
\emd

\subsubsection{The case $a=1$ and $b=1$}
\bbd \label{lemma2.6}
Let $s\ge 0$ and $j=0,\,1$. Then, the following estimates hold for $m \in [1,2)$:
\begin{align}
\big\|\partial_t^j |D|^s \big(K^1_{0\text{\fontshape{n}\selectfont low}}(t,x) \ast u_0(x)\big)(t,\cdot)\big\|_{L^2}&\lesssim (1+t)^{-\frac{n}{2(\sigma-\delta_1)}(\frac{1}{m}-\frac{1}{2})-\frac{s}{2(\sigma-\delta_1)}-j}\|u_0\|_{L^m}, \label{lemma2.6.1} \\
\big\|\partial_t^j |D|^s \big(K^2_{0\text{\fontshape{n}\selectfont low}}(t,x) \ast u_0(x)\big)(t,\cdot)\big\|_{L^2}&\lesssim (1+t)^{-\frac{n}{2\delta_1}(\frac{1}{m}-\frac{1}{2})-\frac{s}{2\delta_1}-j-\frac{\sigma-2\delta_1}{\delta_1}}\|u_0\|_{L^m}, \label{lemma2.6.2}
\end{align}
for any space dimensions $n\ge 1$, and
\begin{align}
\big\|\partial_t^j |D|^s \big(K^1_{1\text{\fontshape{n}\selectfont low}}(t,x) \ast u_1(x)\big)(t,\cdot)\big\|_{L^2}&\lesssim (1+t)^{-\frac{n}{2(\sigma-\delta_1)}(\frac{1}{m}-\frac{1}{2})-\frac{s}{2(\sigma-\delta_1)}-j+\frac{\delta_1}{\sigma-\delta_1}}\|u_1\|_{L^m}, \label{lemma2.6.3} \\
\big\|\partial_t^j |D|^s \big(K^2_{1\text{\fontshape{n}\selectfont low}}(t,x) \ast u_1(x)\big)(t,\cdot)\big\|_{L^2}&\lesssim (1+t)^{-\frac{n}{2\delta_1}(\frac{1}{m}-\frac{1}{2})-\frac{s}{2\delta_1}-j+1}\|u_1\|_{L^m}, \label{lemma2.6.4}
\end{align}
for any space dimensions $n> 2m_0\delta_1$. Moreover, the following estimates hold for $m \in [1,2]$:
\begin{align}
&\big\|\partial_t^j |D|^s \big(K_{0\text{\fontshape{n}\selectfont high}}(t,x) \ast u_0(x)\big)(t,\cdot)\big\|_{L^2} \nonumber \\
&\qquad \lesssim e^{-ct}\Big(t^{-\frac{n}{2(\sigma-\delta_2)}(\frac{1}{m}-\frac{1}{2})-\frac{s-a}{2(\sigma-\delta_2)}-j}+ t^{-\frac{n}{2\delta_2}(\frac{1}{m}-\frac{1}{2})-\frac{s-a}{2\delta_2}-j+ \frac{2\delta_2-\sigma}{\delta_2}}\Big)\|u_0\|_{H^a_m}, \label{lemma2.6.5} \\
&\big\|\partial_t^j |D|^s \big(K_{1\text{\fontshape{n}\selectfont high}}(t,x) \ast u_1(x)\big)(t,\cdot)\big\|_{L^2} \nonumber \\
&\qquad \lesssim e^{-ct}\Big(t^{-\frac{n}{2(\sigma-\delta_2)}(\frac{1}{m}-\frac{1}{2})-\frac{s-a}{2(\sigma-\delta_2)}-j}+ t^{-\frac{n}{2\delta_2}(\frac{1}{m}-\frac{1}{2})-\frac{s-a}{2\delta_2}-j}\Big)\|u_0\|_{H^{[a-2\delta_2]^+}_m}, \label{lemma2.6.6}
\end{align}
for all $a\ge 0$ and for any space dimensions $n\ge 1$, where $c$ is a suitable positive constant.
\ebd
\begin{proof}
Following the proof of Lemma \ref{lemma2.4} we may conclude all the desired estimates in Lemma \ref{lemma2.6} by the aid of the statements in Lemma \ref{lemma2.3}.
\end{proof}

\noindent Hence, we obtain decay properties of Sobolev solutions to (\ref{equation1.1}).
\bmd \label{proposition2.3}
Let $s \ge 0$ and $j=0,1$. Let $m \in [1,2)$. We assume the condition $n>2m_0\delta_1$ and the initial data
$$ (u_0,u_1) \in (L^m \cap H^{s+2j\delta_2}) \times (L^m \cap H^{[s+2(j-1)\delta_2]^+}). $$
Then, the Sobolev solutions to (\ref{equation1.1}) satisfy the following $(L^m \cap L^2)-L^2$ estimates:
\begin{align*}
\big\|\partial_t^j |D|^s u(t,\cdot)\big\|_{L^2} &\lesssim (1+t)^{-\frac{n}{2(\sigma-\delta_1)}(\frac{1}{m}-\frac{1}{2})-\frac{s}{2(\sigma-\delta_1)}-j}\|u_0\|_{L^m \cap H^{s+2j\delta_2}} \\ 
&\qquad + (1+t)^{-\frac{n}{2(\sigma-\delta_1)}(\frac{1}{m}-\frac{1}{2})-\frac{s}{2(\sigma-\delta_1)}-j+ \frac{\delta_1}{\sigma-\delta_1}}\|u_1\|_{L^m \cap H^{[s+2(j-1)\delta_2]^+}}.
\end{align*}
\emd

\section{Proofs} \label{Sec.3}

\subsection{The case $a=1$ and $b=0$} \label{Sec.3.1}
In order to prove Theorem \ref{theorem1.1}, the following auxilliary results come into play.
\bmd \label{proposition3.1}
Let $s\ge 0$ and $j=0,\,1$. Let $m \in [1,2)$. Then, the following estimates hold:
\begin{align}
&\Big\|\partial_t^j |D|^s \Big(\Big(K^1_{0\text{\fontshape{n}\selectfont low}}(t,x)- \mathfrak{F}^{-1}\Big(e^{-t|\xi|^{2(\sigma-\delta_1)}}\chi(\xi)\Big)\Big) \ast u_0(x)\Big)(t,\cdot)\Big\|_{L^2} \nonumber \\
&\qquad \lesssim (1+t)^{-\frac{n}{2(\sigma-\delta_1)}(\frac{1}{m}-\frac{1}{2})-\frac{s}{2(\sigma-\delta_1)}-j-\frac{\sigma-2\delta_1}{\sigma-\delta_1}}\|u_0\|_{L^m}, \label{proposition3.1.1}
\end{align}
for any space dimensions $n\ge 1$, and
\begin{align}
&\Big\|\partial_t^j |D|^s \Big(\Big(K^1_{1\text{\fontshape{n}\selectfont low}}(t,x)- \mathfrak{F}^{-1}\Big(\dfrac{e^{-t|\xi|^{2(\sigma-\delta_1)}}}{|\xi|^{2\delta_1}}\chi(\xi)\Big)\Big) \ast u_1(x)\Big)(t,\cdot)\Big\|_{L^2} \\
&\qquad \lesssim (1+t)^{-\frac{n}{2(\sigma-\delta_1)}(\frac{1}{m}-\frac{1}{2})-\frac{s}{2(\sigma-\delta_1)}-j-\frac{\sigma-3\delta_1}{\sigma-\delta_1}}\|u_1\|_{L^m}, \label{proposition3.1.2}
\end{align}
for any space dimensions $n> 2m_0\delta_1$.
\emd

\begin{proof}
In order to indicate Proposition \ref{proposition3.1}, at first it is reasonable to prove the following estimates:
\begin{align}
&|\xi|^s \chi(|\xi|)\Big|\partial^j_t\Big(\widehat{K^1_0}(t,\xi)- e^{-t|\xi|^{2(\sigma-\delta_1)}}\Big)\Big| \nonumber \\
&\qquad \lesssim e^{-ct|\xi|^{2(\sigma-\delta_1)}}\big(t\,|\xi|^{s+2(2\sigma-3\delta_1)+2j(\sigma-\delta_1)}+ |\xi|^{s+2(\sigma-2\delta_1)+2j(\sigma-\delta_1)}\big), \label{pro3.1.1}
\end{align}
and
\begin{align}
&|\xi|^s \chi(|\xi|)\Big|\partial^j_t \Big(\widehat{K^1_1}(t,\xi)- \frac{e^{-t|\xi|^{2(\sigma-\delta_1)}}}{|\xi|^{2\delta_1}}\Big)\Big| \nonumber \\
&\qquad \lesssim e^{-ct|\xi|^{2(\sigma-\delta_1)}}\big(t\,|\xi|^{s+4(\sigma-2\delta_1)+2j(\sigma-\delta_1)}+ |\xi|^{s+2(\sigma-3\delta_1)+2j(\sigma-\delta_1)}\big), \label{pro3.1.2}
\end{align}
where $c$ is a suitable positive constant. Then, following the proof of Lemma \ref{lemma2.4} we may conclude all the desired statements. In the first step, let us consider (\ref{pro3.1.1}) in the case $j=0$ to present $\widehat{K^1_0}(t,\xi)$ as follows: 
$$ \widehat{K^1_0}(t,\xi)= \frac{-\lambda_2 e^{\lambda_1 t}}{\lambda_1- \lambda_2}= e^{\lambda_1 t}- \frac{\lambda_1 e^{\lambda_1 t}}{\lambda_1- \lambda_2}. $$
By the mean value theorem we obtain
$$ e^{\lambda_1 t}- e^{-t|\xi|^{2(\sigma-\delta_1)}}= t\,\big(\lambda_1+ |\xi|^{2(\sigma-\delta_1)}\big)e^{(\omega\lambda_1-(1-\omega)|\xi|^{2(\sigma-\delta_1)})t}, $$
where $\omega \in [0,1]$. Consequently, we derive
$$ \Big|e^{\lambda_1 t}- e^{-t|\xi|^{2(\sigma-\delta_1)}}\Big| \le t\,\Big|-\lambda_1- |\xi|^{2(\sigma-\delta_1)}\Big| e^{-\min\{-\lambda_1,\,|\xi|^{2(\sigma-\delta_1)}\}t}. $$
Using Newton's binomial theorem we re-write $-\lambda_1$ for small frequencies in the form
\begin{align*}
-\lambda_1 &= \frac{1}{2}|\xi|^{2\delta_1}\Big(1- \sqrt{1-4|\xi|^{2(\sigma-2\delta_1)}}\Big) \\ 
&= \frac{1}{2}|\xi|^{2\delta_1}\Big(1- \Big(1-2|\xi|^{2(\sigma-2\delta_1)}-2|\xi|^{4(\sigma-2\delta_1)}- o\big(|\xi|^{4(\sigma-2\delta_1)}\big)\Big)\Big) \\
&= |\xi|^{2(\sigma-\delta_1)}+ |\xi|^{2(2\sigma-3\delta_1)}+ o\big(|\xi|^{2(2\sigma-3\delta_1)}\big).
\end{align*}
Therefore, we get
\begin{equation} \label{pro3.1.3}
\Big|e^{\lambda_1 t}- e^{-t|\xi|^{2(\sigma-\delta_1)}}\Big| \lesssim  t\,|\xi|^{2(2\sigma-3\delta_1)}e^{-t|\xi|^{2(\sigma-\delta_1)}}.
\end{equation}
Thanks to the asymptotic behavior of characteristic roots for small frequencies, we may arrive at the following estimates:
\begin{equation} \label{pro3.1.4}
|\xi|^s \chi(|\xi|)\Big|\widehat{K^1_0}(t,\xi)- e^{-t|\xi|^{2(\sigma-\delta_1)}}\Big| \lesssim e^{-t|\xi|^{2(\sigma-\delta_1)}}\big(t\,|\xi|^{2(2\sigma-3\delta_1)}+ |\xi|^{2(\sigma-2\delta_1)}\big).
\end{equation}
Hence, the estimate (\ref{pro3.1.1}) is true for $j=0$. Now in oder to estimate (\ref{pro3.1.2}) in the case $j=0$, we can re-write
$$ \widehat{K^1_1}(t,\xi)= \frac{e^{\lambda_1 t}}{\lambda_1- \lambda_2}= \frac{e^{\lambda_1 t}}{|\xi|^{2\delta_1}}- \frac{1}{|\xi|^{2\delta_1}}\Big(1- \frac{|\xi|^{2\delta_1}}{\lambda_1- \lambda_2}\Big)e^{\lambda_1 t}. $$
Hence, we have
$$ \widehat{K^1_1}(t,\xi)- \frac{e^{-t|\xi|^{2(\sigma-\delta_1)}}}{|\xi|^{2\delta_1}}= \frac{e^{\lambda_1 t}- e^{-t|\xi|^{2(\sigma-\delta_1)}}}{|\xi|^{2\delta_1}}- \frac{1}{|\xi|^{2\delta_1}}\Big(1- \frac{|\xi|^{2\delta_1}}{\lambda_1- \lambda_2}\Big)e^{\lambda_1 t}. $$
Moreover, we notice that it holds
$$1- \frac{|\xi|^{2\delta_1}}{\lambda_1- \lambda_2}= 1- \frac{|\xi|^{2\delta_1}}{\sqrt{|\xi|^{4\delta_1}-4|\xi|^{2\sigma}}}= \frac{-4|\xi|^{2\sigma}}{\sqrt{|\xi|^{4\delta_1}-4|\xi|^{2\sigma}}\big(\sqrt{|\xi|^{4\delta_1}-4|\xi|^{2\sigma}}+ |\xi|^{2\delta_1}\big)}\sim -|\xi|^{2(\sigma- 2\delta_1)}. $$
Using again the estimate (\ref{pro3.1.3}) and the asymptotic behavior of characteristic roots for small frequencies we may conclude
\begin{equation} \label{pro3.1.5}
|\xi|^s \chi(|\xi|)\Big|\widehat{K^1_1}(t,\xi)- \frac{e^{-t|\xi|^{2(\sigma-\delta_1)}}}{|\xi|^{2\delta_1}}\Big| \lesssim e^{-ct|\xi|^{2(\sigma-\delta_1)}}\big(t\,|\xi|^{s+4(\sigma-2\delta_1)}+ |\xi|^{s+2(\sigma-3\delta_1)}\big).
\end{equation}
Therefore, the estimate (\ref{pro3.1.2}) is true for $j=0$. By analogous arguments we also obtain the following estimates for $j=1$:
\begin{align}
|\xi|^s \chi(|\xi|)\Big|\partial_t\Big(\widehat{K^1_0}(t,\xi)- e^{-t|\xi|^{2(\sigma-\delta_1)}}\Big)\Big| &\lesssim  e^{-t|\xi|^{2(\sigma-\delta_1)}}\big(t\,|\xi|^{2(3\sigma-4\delta_1)}+ |\xi|^{2(2\sigma-3\delta_1)}\big), \label{pro3.1.6} \\ 
|\xi|^s \chi(|\xi|)\Big|\partial_t \Big(\widehat{K^1_1}(t,\xi)- \frac{e^{-t|\xi|^{2(\sigma-\delta_1)}}}{|\xi|^{2\delta_1}}\Big)\Big| &\lesssim e^{-ct|\xi|^{2(\sigma-\delta_1)}}\big(t\,|\xi|^{s+2(4\sigma-5\delta_1)}+ |\xi|^{s+4(\sigma-2\delta_1)}\big). \label{pro3.1.7}
\end{align}
Thus, it is obvious that all the estimates from (\ref{pro3.1.4}) to (\ref{pro3.1.7}) imply immediately (\ref{pro3.1.1}) and (\ref{pro3.1.2}). This completes our proof.
\end{proof}

\bmd \label{proposition3.2}
Let $s\ge 0$ and $j=0,\,1$. Let $m \in [1,2)$. Then, the following estimates hold:
\begin{align}
&\Big\|\partial_t^j |D|^s \Big(\Big(K^1_0(t,x)- \mathfrak{F}^{-1}\Big(e^{-t|\xi|^{2(\sigma-\delta_1)}}\Big)\Big) \ast u_0(x)\Big)(t,\cdot)\Big\|_{L^2} \nonumber \\
&\qquad \lesssim (1+t)^{-\frac{n}{2(\sigma-\delta_1)}(\frac{1}{m}-\frac{1}{2})-\frac{s}{2(\sigma-\delta_1)}-j-\frac{\sigma-2\delta_1}{\sigma-\delta_1}}\|u_0\|_{L^m}+ e^{-ct}t^{-\frac{n}{2\delta_1}(\frac{1}{m_1}-\frac{1}{2})-\frac{s+j\sigma-a_1}{2\delta_1}}\|u_0\|_{H^{a_1}_{m_1}} \nonumber \\
&\hspace{8cm} +e^{-ct}t^{-\frac{n}{2(\sigma-\delta_1)}(\frac{1}{m_2}-\frac{1}{2})-\frac{s-a_2}{2(\sigma-\delta_1)}-j}\|u_0\|_{H^{a_2}_{m_2}} \label{proposition3.2.1}
\end{align}
for any space dimensions $n\ge 1$, and
\small
\begin{align}
&\Big\|\partial_t^j |D|^s \Big(\Big(K^1_1(t,x)- \mathfrak{F}^{-1}\Big(\dfrac{e^{-t|\xi|^{2(\sigma-\delta_1)}}}{|\xi|^{2\delta_1}}\Big)\Big) \ast u_1(x)\Big)(t,\cdot)\Big\|_{L^2} \nonumber \\
&\qquad \lesssim (1+t)^{-\frac{n}{2(\sigma-\delta_1)}(\frac{1}{m}-\frac{1}{2})-\frac{s}{2(\sigma-\delta_1)}-j-\frac{\sigma-3\delta_1}{\sigma-\delta_1}}\|u_1\|_{L^m}+ e^{-ct}t^{-\frac{n}{2\delta_1}(\frac{1}{m_1}-\frac{1}{2})-\frac{s+j\sigma-a_1}{2\delta_1}}\|u_1\|_{H^{[a_1-\sigma]^+}_{m_1}} \nonumber \\
&\hspace{8cm} +e^{-ct}t^{-\frac{n}{2(\sigma-\delta_1)}(\frac{1}{m_2}-\frac{1}{2})-\frac{s-a_2}{2(\sigma-\delta_1)}-j}\|u_1\|_{H^{[a_2-2\delta_1]^+}_{m_2}} \label{proposition3.2.2}
\end{align}
\normalsize
for any space dimensions $n> 2m_0\delta_1$. Here $a_1,\,a_2 \ge 0$, $m_1,\,m_2 \in [1,2]$ and $c$ is a suitable positive constant.
\emd

\begin{proof}
At first, let us re-write the expression in the $L^2$ norm of (\ref{proposition3.2.1}) as follows:
\begin{align*}
&\partial_t^j |D|^s \Big(\Big(K^1_0(t,x)- \mathfrak{F}^{-1}\Big(e^{-t|\xi|^{2(\sigma-\delta_1)}}\Big)\Big) \ast u_0(x)\Big) \\ 
&\qquad= \partial_t^j |D|^s \Big(\Big(K^1_{0\text{\fontshape{n}\selectfont low}}(t,x)- \mathfrak{F}^{-1}\Big(e^{-t|\xi|^{2(\sigma-\delta_1)}}\chi(\xi)\Big)\Big) \ast u_0(x)\Big) \\
&\qquad \quad+ \partial_t^j |D|^s \big(K^1_{0\text{\fontshape{n}\selectfont high}}(t,x) \ast u_0(x)\big)+ \partial_t^j |D|^s \Big(\mathfrak{F}^{-1}\Big(e^{-t|\xi|^{2(\sigma-\delta_1)}}\big(1-\chi(\xi)\big)\Big) \ast u_0(x)\Big).
\end{align*}
We notice that it holds
$$ |\xi|^s \big(1-\chi(\xi)\big) \big|\partial_t^j e^{-t|\xi|^{2(\sigma-\delta_1)}}\big|\lesssim |\xi|^{s+2j(\sigma-\delta_1)}e^{-t|\xi|^{2(\sigma-\delta_1)}}. $$
Following the proof of Lemma \ref{lemma2.4} we derive
\begin{equation}
\Big\|\partial_t^j |D|^s \Big(\mathfrak{F}^{-1}\Big(e^{-t|\xi|^{2(\sigma-\delta_1)}}\big(1-\chi(\xi)\big)\Big) \ast u_0(x)\Big)(t,\cdot)\Big\|_{L^2}\lesssim e^{-ct}t^{-\frac{n}{2(\sigma-\delta_1)}(\frac{1}{m_2}-\frac{1}{2})-j-\frac{s-a_2}{2(\sigma-\delta_1)}}\|u_0\|_{H^{a_2}_{m_2}}. \label{pro3.2.1}
\end{equation}
Therefore, combining (\ref{lemma2.4.5}), (\ref{proposition3.1.1}) and (\ref{pro3.2.1}) we may arrive at (\ref{proposition3.2.1}). In an analogous way to get (\ref{pro3.2.1}) we also obtain
\small
\begin{equation}
\Big\|\partial_t^j |D|^s \Big(\mathfrak{F}^{-1}\Big(\dfrac{e^{-t|\xi|^{2(\sigma-\delta_1)}}}{|\xi|^{2\delta_1}}\big(1-\chi(\xi)\big)\Big) \ast u_1(x)\Big)(t,\cdot)\Big\|_{L^2}\lesssim e^{-ct}t^{-\frac{n}{2(\sigma-\delta_1)}(\frac{1}{m_2}-\frac{1}{2})-\frac{s-a_2}{2(\sigma-\delta_1)}-j}\|u_1\|_{H^{[a_2-2\delta_1]^+}_{m_2}}. \label{pro3.2.2}
\end{equation}
\normalsize
Hence, combining (\ref{lemma2.4.6}), (\ref{proposition3.1.2}) and (\ref{pro3.2.2}) we may conclude (\ref{proposition3.2.2}). Summurizing, the proof of Proposition \ref{proposition3.2} is completed.
\end{proof}

\begin{proof}[Proof of Theorem \ref{theorem1.1}]
In order to prove the asymptotic profile of solutions to (\ref{equation1.1}), we may estimate
\begin{align*}
&\Big\|\partial_t^j |D|^s \Big(u(t,\cdot)- P_1\,\mathfrak{F}^{-1}\Big(\dfrac{e^{-t|\xi|^{2(\sigma-\delta_1)}}}{|\xi|^{2\delta_1}}\Big)(t,\cdot)\Big)(t,\cdot)\Big\|_{L^2} \\
&\qquad \lesssim \big\|\partial_t^j |D|^s \big(K^1_0(t,x) \ast u_0(x)\big)(t,\cdot)\big\|_{L^2}+ \big\|\partial_t^j |D|^s \big(K^2_0(t,x) \ast u_0(x)\big)(t,\cdot)\big\|_{L^2} \\
&\qquad \quad+ \big\||D|^s \big(K^2_1(t,x) \ast u_1(x)\big)(t,\cdot)\big\|_{L^2}+ \Big\|\partial_t^j|D|^s \Big(\Big(K^1_1(t,x)- \mathfrak{F}^{-1}\Big(\frac{e^{-t|\xi|^{2(\sigma-\delta_1)}}}{|\xi|^{2\delta_1}}\Big)\Big) \ast u_1(x)\Big)(t,\cdot)\Big\|_{L^2} \\
&\qquad \quad+ \Big\|\partial_t^j |D|^s \Big(\mathfrak{F}^{-1}\Big(\frac{e^{-t|\xi|^{2(\sigma-\delta_1)}}}{|\xi|^{2\delta_1}}\Big) \ast u_1(x)- P_1\,\mathfrak{F}^{-1}\Big(\dfrac{e^{-t|\xi|^{2(\sigma-\delta_1)}}}{|\xi|^{2\delta_1}}\Big)\Big)(t,\cdot)\Big\|_{L^2} \\
&\qquad =: I_1+ I_2+ I_3+ I_4+ I_5.
\end{align*}
Combining (\ref{lemma2.4.1}) and (\ref{lemma2.4.5}), (\ref{lemma2.4.2}) and (\ref{lemma2.4.5}), (\ref{lemma2.4.4}) and (\ref{lemma2.4.6}) we get
\begin{align*}
I_1 &\lesssim (1+t)^{-\frac{n}{4(\sigma-\delta_1)}- \frac{s}{2(\sigma-\delta_1)}-j} \|u_0\|_{\mathcal{A}^1_0}, \\ 
I_2 &\lesssim (1+t)^{-\frac{n}{4\delta_1}-\frac{s}{2\delta_1}-j- \frac{\sigma-2\delta_1}{\delta_1}} \|u_0\|_{\mathcal{A}^1_0}, \\
I_3 &\lesssim (1+t)^{-\frac{n}{4\delta_1}- \frac{s}{2\delta_1}-j+1} \|u_1\|_{\mathcal{A}^1_1},
\end{align*}
respectively. By (\ref{proposition3.2.2}) we obtain
$$ I_4 \lesssim (1+t)^{-\frac{n}{2(\sigma-\delta_1)}(\frac{1}{m}-\frac{1}{2})-\frac{s}{2(\sigma-\delta_1)}-j-\frac{\sigma-3\delta_1}{\sigma-\delta_1}}\|u_1\|_{\mathcal{A}^1_1}. $$
To control $I_5$, we shall apply Lemma \ref{L^1.Lemma}. Indeed, at first it is clear that using Parseval-Plancherel formula and the change of variables $|\xi|= t^{-\frac{1}{2(\sigma-\delta_1)}}|\eta|$ gives
\begin{align}
\Big\|\partial^j_t |D|^s \mathfrak{F}^{-1}\Big(\dfrac{e^{-t|\xi|^{2(\sigma-\delta_1)}}}{|\xi|^{2\delta_1}}\Big)(t,\cdot)\Big\|_{L^2} &= \Big\||\xi|^{s+2j(\sigma-\delta_1)-2\delta_1}e^{-t|\xi|^{2(\sigma-\delta_1)}}\Big\|_{L^2} \nonumber \\
&= t^{-\frac{n}{4(\sigma-\delta_1)}-\frac{s}{2(\sigma-\delta_1)}-j+\frac{\delta_1}{\sigma-\delta_1}} \Big(\int_0^\ity |\eta|^{2s-4\delta_1+4j(\sigma-\delta_1)}e^{-2|\eta|^{2(\sigma-\delta_1)}d\eta}\Big)^{\frac{1}{2}} \nonumber \\
&= C\,t^{-\frac{n}{4(\sigma-\delta_1)}-\frac{s}{2(\sigma-\delta_1)}-j+\frac{\delta_1}{\sigma-\delta_1}} \label{the1.1.1}
\end{align}
with the constant $C:= \big(\int_0^\ity |\eta|^{2s-4\delta_1+4j(\sigma-\delta_1)}e^{-2|\eta|^{2(\sigma-\delta_1)}d\eta}\big)^{\frac{1}{2}}> 0$, where we used the condition $n> 4\delta_1$. Then, we employ Lemma \ref{L^1.Lemma} to derive
$$ I_5= o\Big(t^{-\frac{n}{4(\sigma-\delta_1)}-\frac{s}{2(\sigma-\delta_1)}-j+ \frac{\delta_1}{\sigma-\delta_1}}\Big) \quad \text{ as }t\to \ity. $$
Therefore, from all the above estimates for $I_k$ with $k=1,\cdots,5$ and paying attention that it holds $ -\frac{n}{4\delta_1}- \frac{s}{2\delta_1}+1< -\frac{n}{4(\sigma-\delta_1)}-\frac{s}{2(\sigma-\delta_1)}+\frac{\delta_1}{\sigma-\delta_1}$
with the condtion $n> 4\delta_1$ we may conclude immediately (\ref{theorem1.1.1}). Then, from (\ref{theorem1.1.1}) and (\ref{the1.1.1}) we may arrive at the desired upper bound in the following way:
\begin{align*}
\big\|\partial^j_t |D|^s u(t,\cdot)\big\|_{L^2} &\le \Big\|\partial_t^j |D|^s \Big(u(t,\cdot)- P_1\,\mathfrak{F}^{-1}\Big(\dfrac{e^{-t|\xi|^{2(\sigma-\delta_1)}}}{|\xi|^{2\delta_1}}\Big)(t,\cdot)\Big)\Big\|_{L^2} \\
&\qquad+ |P_1|\,\Big\|\partial^j_t |D|^s \mathfrak{F}^{-1}\Big(\dfrac{e^{-t|\xi|^{2(\sigma-\delta_1)}}}{|\xi|^{2\delta_1}}\Big)(t,\cdot)\Big\|_{L^2} \\
&\le C|P_1|\,t^{-\frac{n}{4(\sigma-\delta_1)}-\frac{s}{2(\sigma-\delta_1)}-j+\frac{\delta_1}{\sigma-\delta_1}}+ o\Big(t^{-\frac{n}{4(\sigma-\delta_1)}-\frac{s}{2(\sigma-\delta_1)}-j+\frac{\delta_1}{\sigma-\delta_1}}\Big) \\
&\le C_2\,t^{-\frac{n}{4(\sigma-\delta_1)}-\frac{s}{2(\sigma-\delta_1)}-j+\frac{\delta_1}{\sigma-\delta_1}} \quad \text{ as }t\to \ity, 
\end{align*}
where $C_2$ is a suitable positive constant. Moreover, to indicate the lower bound, we can proceed as follows:
\begin{align*}
\big\|\partial^j_t |D|^s u(t,\cdot)\big\|_{L^2}&\ge |P_1|\,\Big\|\partial^j_t |D|^s \mathfrak{F}^{-1}\Big(\dfrac{e^{-t|\xi|^{2(\sigma-\delta_1)}}}{|\xi|^{2\delta_1}}\Big)(t,\cdot)\Big\|_{L^2} \\
&\qquad- \Big\|\partial_t^j |D|^s \Big(u(t,\cdot)- P_1\,\mathfrak{F}^{-1}\Big(\dfrac{e^{-t|\xi|^{2(\sigma-\delta_1)}}}{|\xi|^{2\delta_1}}\Big)(t,\cdot)\Big)\Big\|_{L^2} \\ 
&\ge C|P_1|\,t^{-\frac{n}{4(\sigma-\delta_1)}-\frac{s}{2(\sigma-\delta_1)}-j+ \frac{\delta_1}{\sigma-\delta_1}}- o\Big(t^{-\frac{n}{4(\sigma-\delta_1)}-j-\frac{s}{2(\sigma-\delta_1)}+ \frac{\delta_1}{\sigma-\delta_1}}\Big) \\
&\ge C_1\,t^{-\frac{n}{4(\sigma-\delta_1)}-\frac{s}{2(\sigma-\delta_1)}-j+ \frac{\delta_1}{\sigma-\delta_1}}\quad \text{ as }t\to \ity,
\end{align*}
where $C_1$ is a suitable positive constant. Summurizing, Theorem \ref{theorem1.1} is proved.
\end{proof}

\subsection{The case $a=0$ and $b=1$} \label{Sec.3.2}
In order to prove Theorem \ref{theorem1.2}, the following auxilliary results come into play.
\bmd \label{proposition3.3}
Let $s\ge 0$ and $m \in [1,2)$. Then, the following estimates hold:
\begin{align}
&\Big\||D|^s \Big(\Big(K^{\cos}_{0\text{\fontshape{n}\selectfont low}}(t,x)- \mathfrak{F}^{-1}\Big(e^{-t|\xi|^{2\delta_2}}\cos\big(t|\xi|^\sigma\big)\chi(\xi)\Big)\Big) \ast u_0(x)\Big)(t,\cdot)\Big\|_{L^2} \nonumber \\
&\qquad \lesssim (1+t)^{-\frac{n}{2\delta_2}(\frac{1}{m}-\frac{1}{2})-1-\frac{s}{2\delta_2}+\frac{\sigma}{2\delta_2}}\|u_0\|_{L^m}, \label{proposition3.3.1}
\end{align}
\begin{align}
&\Big\||D|^s \Big(\Big(\partial_t K^{\cos}_{0\text{\fontshape{n}\selectfont low}}(t,x)+ \mathfrak{F}^{-1}\Big(|\xi|^\sigma e^{-t|\xi|^{2\delta}}\sin\big(t|\xi|^\sigma\big)\chi(\xi)\Big)\Big) \ast u_0(x)\Big)(t,\cdot)\Big\|_{L^2} \nonumber \\
&\qquad \lesssim (1+t)^{-\frac{n}{2\delta_2}(\frac{1}{m}-\frac{1}{2})-1-\frac{s}{2\delta_2}}\|u_0\|_{L^m}, \label{proposition3.3.2} \\
&\Big\||D|^s \Big(\Big(\partial_t K_{1\text{\fontshape{n}\selectfont low}}(t,x)- \mathfrak{F}^{-1}\Big(e^{-t|\xi|^{2\delta_2}}\cos\big(t|\xi|^\sigma\big)\chi(\xi)\Big)\Big) \ast u_1(x)\Big)(t,\cdot)\Big\|_{L^2} \nonumber \\
&\qquad \lesssim (1+t)^{-\frac{n}{2\delta_2}(\frac{1}{m}-\frac{1}{2})-1-\frac{s}{2\delta_2}+\frac{\sigma}{2\delta_2}}\|u_1\|_{L^m}, \label{proposition3.3.3}
\end{align}
for any space dimensions $n\ge 1$, and
\begin{align}
&\Big\||D|^s \Big(\Big(K_{1\text{\fontshape{n}\selectfont low}}(t,x)- \mathfrak{F}^{-1}\Big(e^{-t|\xi|^{2\delta_2}}\dfrac{\sin\big(t|\xi|^\sigma\big)}{|\xi|^\sigma}\chi(\xi)\Big)\Big) \ast u_1(x)\Big)(t,\cdot)\Big\|_{L^2} \nonumber \\
&\qquad \lesssim (1+t)^{-\frac{n}{2\delta_2}(\frac{1}{m}-\frac{1}{2})-1-\frac{s}{2\delta_2}+\frac{\sigma}{\delta_2}}\|u_1\|_{L^m}, \label{proposition3.3.4}
\end{align}
for any space dimensions $n> m_0\sigma$.
\emd

\begin{proof}
To prove Proposition \ref{proposition3.3}, we only need to show the following estimates:
\begin{align}
|\xi|^s \chi(|\xi|)\Big|\widehat{K^{\cos}_0}(t,\xi)- e^{-\frac{1}{2}t|\xi|^{2\delta_2}}\cos(t|\xi|^\sigma)\Big| &\lesssim t\,e^{-ct|\xi|^{2\delta_2}}|\xi|^{s+4\delta_2-\sigma}, \label{pro3.3.1} \\
|\xi|^s \chi(|\xi|)\Big|\widehat{K_1}(t,\xi)- e^{-\frac{1}{2}t|\xi|^{2\delta_2}}\dfrac{\sin(t|\xi|^\sigma)}{|\xi|^\sigma}\Big| &\lesssim e^{-ct|\xi|^{2\delta_2}}\big(t\,|\xi|^{s+4\delta_2-2\sigma}+|\xi|^{s+4\delta_2-3\sigma}\big), \label{pro3.3.2}
\end{align}
and
\begin{align}
|\xi|^s \chi(|\xi|)\Big|\partial_t \widehat{K^{\cos}_0}(t,\xi)+ e^{-\frac{1}{2}t|\xi|^{2\delta_2}}|\xi|^\sigma \sin(t|\xi|^\sigma)\Big| &\lesssim e^{-ct|\xi|^{2\delta_2}}\big(t\,|\xi|^{s+4\delta_2}+ |\xi|^{s+2\delta_2}\big), \label{pro3.3.3} \\
|\xi|^s \chi(|\xi|)\Big|\partial_t \widehat{K_1}(t,\xi)- e^{-\frac{1}{2}t|\xi|^{2\delta_2}}\cos(t|\xi|^\sigma)\Big| &\lesssim e^{-ct|\xi|^{2\delta_2}}\big(t\,|\xi|^{s+4\delta_2-\sigma}+|\xi|^{s+2\delta_2-\sigma}\big), \label{pro3.3.4}
\end{align}
where $c$ is a suitable positive constant. Then, following the proof of Lemma \ref{lemma2.4} we may conclude all the desired statements. Indeed, at first let us consider (\ref{pro3.3.1}) to re-wirte
\begin{equation}
\widehat{K^{\cos}_0}(t,\xi)- e^{-\frac{1}{2}t|\xi|^{2\delta_2}}\cos(t|\xi|^\sigma)= e^{-\frac{1}{2}t|\xi|^{2\delta_2}}\big(\cos\big(t|\xi|^\sigma f(|\xi|)\big)- \cos(t|\xi|^\sigma)\big). \label{pro3.3.5}
\end{equation}
Using the mean value theorem we have
\begin{equation}
\cos\big(t|\xi|^\sigma f(|\xi|)\big)- \cos(t|\xi|^\sigma)= -t|\xi|^\sigma \big(f(|\xi|)-1\big) \sin\big(t|\xi|^\sigma \big(\omega_1 f(|\xi|)+ 1-\omega_1\big)\big), \label{pro3.3.6}
\end{equation}
where $\omega_1 \in [0,1]$. Moreover, it holds for small frequencies
\begin{equation}
f(|\xi|)-1= -\frac{|\xi|^{4\delta_2-2\sigma}}{4\big(1+ \sqrt{1-\frac{1}{4}|\xi|^{4\delta_2-2\sigma}}\big)}. \label{pro3.3.7}
\end{equation}
From (\ref{pro3.3.5}) to (\ref{pro3.3.7}), we imply immediately (\ref{pro3.3.1}). In order to indicate (\ref{pro3.3.2}), we get
\small
\begin{align}
\widehat{K_1}(t,\xi)- e^{-\frac{1}{2}t|\xi|^{2\delta_2}}\dfrac{\sin(t|\xi|^\sigma)}{|\xi|^\sigma}&= e^{-\frac{1}{2}t|\xi|^{2\delta_2}}\Big(\dfrac{\sin\big(t|\xi|^\sigma f(|\xi|)\big)}{|\xi|^\sigma f(|\xi|)}- \dfrac{\sin(t|\xi|^\sigma)}{|\xi|^\sigma}\Big) \nonumber \\ 
&= \dfrac{e^{-\frac{1}{2}t|\xi|^{2\delta_2}}}{|\xi|^\sigma}\Big(\dfrac{\sin\big(t|\xi|^\sigma f(|\xi|)\big)- \sin(t|\xi|^\sigma)}{f(|\xi|)}- \sin(t|\xi|^\sigma)\dfrac{f(|\xi|)-1}{f(|\xi|)}\Big). \label{pro3.3.8}
\end{align}
\normalsize
Applying again the mean value theorem leads to
\begin{equation}
\sin\big(t|\xi|^\sigma f(|\xi|)\big)- \sin(t|\xi|^\sigma)= t|\xi|^\sigma \big(f(|\xi|)-1\big) \cos\big(t|\xi|^\sigma \big(\omega_2 f(|\xi|)+ 1-\omega_2\big)\big), \label{pro3.3.9}
\end{equation}
where $\omega_2 \in [0,1]$. Hence, combining from (\ref{pro3.3.7}) to (\ref{pro3.3.9}) gives (\ref{pro3.3.2}). Next, to show (\ref{pro3.3.3}), we notice that it holds
\begin{align}
&\partial_t \widehat{K^{\cos}_0}(t,\xi)+ e^{-\frac{1}{2}t|\xi|^{2\delta_2}}|\xi|^\sigma \sin(t|\xi|^\sigma) \nonumber \\ 
&\qquad = -e^{-\frac{1}{2}t|\xi|^{2\delta_2}}\Big(\frac{1}{2}|\xi|^{2\delta_2} \cos\big(t|\xi|^\sigma f(|\xi|)\big)+ |\xi|^\sigma f(|\xi|)\sin\big(t|\xi|^\sigma f(|\xi|)\big)- |\xi|^\sigma \sin(t|\xi|^\sigma)\Big) \nonumber \\
&\qquad = -e^{-\frac{1}{2}t|\xi|^{2\delta_2}}\Big(\frac{1}{2}|\xi|^{2\delta_2} \cos\big(t|\xi|^\sigma f(|\xi|)\big)+ |\xi|^\sigma f(|\xi|)\big(\sin\big(t|\xi|^\sigma f(|\xi|)\big)- \sin(t|\xi|^\sigma)\big) \nonumber \\
&\hspace{9cm}+ |\xi|^\sigma\big(f(|\xi|)-1\big)\sin(t|\xi|^\sigma) \Big). \label{pro3.3.10}
\end{align}
From (\ref{pro3.3.7}), (\ref{pro3.3.9}) and (\ref{pro3.3.10}), we may arrive at (\ref{pro3.3.3}). Finally, by using the relation
\begin{align*}
&\partial_t \widehat{K_1}(t,\xi)- e^{-\frac{1}{2}t|\xi|^{2\delta_2}}\cos(t|\xi|^\sigma) \\
&\qquad= e^{-\frac{1}{2}t|\xi|^{2\delta_2}}\Big(-|\xi|^{2\delta_2- \sigma} \dfrac{\sin\big(t|\xi|^\sigma f(|\xi|)\big)}{f(|\xi|)}+ \big(\cos\big(t|\xi|^\sigma f(|\xi|)\big)- \cos(t|\xi|^\sigma)\big)\Big)
\end{align*}
and combining (\ref{pro3.3.6}), (\ref{pro3.3.7}) we may conclude (\ref{pro3.3.4}). This completes the proof of Proposition \ref{proposition3.3}.
\end{proof}

\bmd \label{proposition3.4}
Let $s\ge 0$ and Let $m \in [1,2)$. Then, the following estimates hold:
\begin{align}
&\Big\||D|^s \Big(\Big(K^{\cos}_0(t,x)- \mathfrak{F}^{-1}\Big(e^{-t|\xi|^{2\delta_2}}\cos\big(t|\xi|^\sigma\big)\Big)\Big) \ast u_0(x)\Big)(t,\cdot)\Big\|_{L^2} \nonumber \\
&\qquad \lesssim (1+t)^{-\frac{n}{2\delta_2}(\frac{1}{m}-\frac{1}{2})-1-\frac{s}{2\delta_2}+\frac{\sigma}{2\delta_2}}\|u_0\|_{L^m}+ e^{-ct}t^{-\frac{n}{2(\sigma-\delta_2)}(\frac{1}{m_1}-\frac{1}{2})-\frac{s-a_1}{2(\sigma-\delta_2)}}\|u_0\|_{H^{a_1}_{m_1}} \nonumber \\
&\qquad \quad +e^{-ct}t^{-\frac{n}{2\delta_2}(\frac{1}{m_2}-\frac{1}{2})-\frac{s-a_2}{2\delta_2}+\frac{2\delta_2-\sigma}{\delta_2}}\|u_0\|_{H^{a_2}_{m_2}}+ e^{-ct}t^{-\frac{n}{2\delta_2}(\frac{1}{m_3}-\frac{1}{2})-\frac{s-a_3}{2\delta_2}}\|u_0\|_{H^{a_3}_{m_3}}, \label{proposition3.4.1} \\
&\Big\||D|^s \Big(\Big(\partial_t K^{\cos}_0(t,x)+ \mathfrak{F}^{-1}\Big(|\xi|^\sigma e^{-t|\xi|^{2\delta_2}}\sin\big(t|\xi|^\sigma\big)\Big)\Big) \ast u_0(x)\Big)(t,\cdot)\Big\|_{L^2} \nonumber \\
&\qquad \lesssim (1+t)^{-\frac{n}{2\delta_2}(\frac{1}{m}-\frac{1}{2})-1-\frac{s}{2\delta_2}}\|u_0\|_{L^m}+ e^{-ct}t^{-\frac{n}{2(\sigma-\delta_2)}(\frac{1}{m_1}-\frac{1}{2})-1-\frac{s-a_1}{2(\sigma-\delta_2)}}\|u_0\|_{H^{a_1}_{m_1}} \nonumber \\
&\qquad \quad + e^{-ct}t^{-\frac{n}{2\delta_2}(\frac{1}{m_2}-\frac{1}{2})-1-\frac{s-a_2}{2\delta_2}+\frac{2\delta_2-\sigma}{\delta_2}}\|u_0\|_{H^{a_2}_{m_2}}+ e^{-ct}t^{-\frac{n}{2\delta_2}(\frac{1}{m_3}-\frac{1}{2})-\frac{s-a_3}{2\delta_2}-\frac{\sigma}{2\delta_2}}\|u_0\|_{H^{a_3}_{m_3}}, \label{proposition3.4.2} \\
&\Big\||D|^s \Big(\Big(\partial_t K_1(t,x)- \mathfrak{F}^{-1}\Big(e^{-t|\xi|^{2\delta_2}}\cos\big(t|\xi|^\sigma\big)\Big)\Big) \ast u_1(x)\Big)(t,\cdot)\Big\|_{L^2} \nonumber \\
&\qquad \lesssim (1+t)^{-\frac{n}{2\delta_2}(\frac{1}{m}-\frac{1}{2})-1-\frac{s}{2\delta_2}+\frac{\sigma}{2\delta_2}}\|u_1\|_{L^m}+ e^{-ct}t^{-\frac{n}{2(\sigma-\delta_2)}(\frac{1}{m_1}-\frac{1}{2})-1-\frac{s-a_1}{2(\sigma-\delta_2)}}\|u_1\|_{H^{[a_1-2\delta_2]^+}_{m_1}} \nonumber \\
&\qquad \quad + e^{-ct}t^{-\frac{n}{2\delta_2}(\frac{1}{m_2}-\frac{1}{2})-1-\frac{s-a_2}{2\delta_2}}\|u_1\|_{H^{[a_2-2\delta_2]^+}_{m_2}}+ e^{-ct}t^{-\frac{n}{2\delta_2}(\frac{1}{m_3}-\frac{1}{2})-\frac{s-a_3}{2\delta_2}}\|u_1\|_{H^{a_3}_{m_3}}, \label{proposition3.4.3}
\end{align}
for any space dimensions $n\ge 1$, and
\begin{align}
&\Big\||D|^s \Big(\Big(K_1(t,x)- \mathfrak{F}^{-1}\Big(e^{-t|\xi|^{2\delta_2}}\dfrac{\sin\big(t|\xi|^\sigma\big)}{|\xi|^\sigma}\Big)\Big) \ast u_1(x)\Big)(t,\cdot)\Big\|_{L^2} \nonumber \\
&\qquad \lesssim (1+t)^{-\frac{n}{2\delta_2}(\frac{1}{m}-\frac{1}{2})-1-\frac{s}{2\delta_2}+\frac{\sigma}{\delta_2}}\|u_1\|_{L^m}+ e^{-ct}t^{-\frac{n}{2(\sigma-\delta_2)}(\frac{1}{m_1}-\frac{1}{2})-\frac{s-a_1}{2(\sigma-\delta_2)}}\|u_1\|_{H^{[a_1-2\delta_2]^+}_{m_1}} \nonumber \\
&\qquad \quad + e^{-ct}t^{-\frac{n}{2\delta_2}(\frac{1}{m_2}-\frac{1}{2})-\frac{s-a_2}{2\delta_2}}\|u_1\|_{H^{[a_2-2\delta_2]^+}_{m_2}}+ e^{-ct}t^{-\frac{n}{2\delta_2}(\frac{1}{m_3}-\frac{1}{2})-\frac{s-a_3}{2\delta_2}}\|u_1\|_{H^{[a_3-\sigma]^+}_{m_3}}, \label{proposition3.4.4}
\end{align}
for any space dimensions $n> m_0\sigma$. Here $a_1,\,a_2,\,a_3 \ge 0$, $m_1,\,m_2,\,m_3 \in [1,2]$ and $c$ is a suitable positive constant.
\emd

\begin{proof}
At first, let us re-write the expression in the $L^2$ norm of (\ref{proposition3.4.1}) as follows:
\begin{align*}
&|D|^s \Big(\Big(K^{\cos}_0(t,x)- \mathfrak{F}^{-1}\Big(e^{-t|\xi|^{2\delta_2}}\cos\big(t|\xi|^\sigma\big)\Big)\Big) \ast u_0(x)\Big) \\ 
&\qquad= |D|^s \Big(\Big(K^{\cos}_{0\text{\fontshape{n}\selectfont low}}(t,x)- \mathfrak{F}^{-1}\Big(e^{-t|\xi|^{2\delta_2}}\cos\big(t|\xi|^\sigma\big)\chi(\xi)\Big)\Big) \ast u_0(x)\Big) \\
&\qquad \quad+ |D|^s \big(K_{0\text{\fontshape{n}\selectfont high}}(t,x) \ast u_0(x)\big)+ |D|^s \Big(\mathfrak{F}^{-1}\Big(e^{-t|\xi|^{2\delta_2}}\cos\big(t|\xi|^\sigma\big)\big(1-\chi(\xi)\big)\Big) \ast u_0(x)\Big)
\end{align*}
It holds that
$$ |\xi|^s \big(1-\chi(\xi)\big) \big|e^{-t|\xi|^{2\delta_2}}\cos\big(t|\xi|^\sigma\big)\big|\lesssim |\xi|^s e^{-t|\xi|^{2\delta_2}}. $$
Following the proof of Lemma \ref{lemma2.4} we obtain
\begin{equation}
\Big\||D|^s \Big(\mathfrak{F}^{-1}\Big(e^{-t|\xi|^{2\delta_2}}\cos\big(t|\xi|^\sigma\big)\big(1-\chi(\xi)\big)\Big) \ast u_0(x)\Big)(t,\cdot)\Big\|_{L^2}\lesssim e^{-ct}t^{-\frac{n}{2\delta_2}(\frac{1}{m_3}-\frac{1}{2})-\frac{s-a_3}{2\delta_2}}\|u_0\|_{H^{a_3}_{m_3}}. \label{pro3.4.1}
\end{equation}
Therefore, combining (\ref{lemma2.5.4}), (\ref{proposition3.3.1}) and (\ref{pro3.4.1}) we may arrive at (\ref{proposition3.4.1}). In analogous ways to get (\ref{pro3.2.1}) we derive the following estimates:
\begin{align}
\Big\||D|^s \Big(\mathfrak{F}^{-1}\Big(e^{-t|\xi|^{2\delta_2}}\sin\big(t|\xi|^\sigma\big)\big(1-\chi(\xi)\big)\Big) \ast u_0(x)\Big)(t,\cdot)\Big\|_{L^2} &\lesssim e^{-ct}t^{-\frac{n}{2\delta_2}(\frac{1}{m_3}-\frac{1}{2})-\frac{s-a_3}{2\delta_2}}\|u_{01}\|_{H^{a_3}_{m_3}}, \label{pro3.4.2} \\
\Big\||D|^s \Big(\mathfrak{F}^{-1}\Big(e^{-t|\xi|^{2\delta_2}}\cos\big(t|\xi|^\sigma\big)\big(1-\chi(\xi)\big)\Big) \ast u_1(x)\Big)(t,\cdot)\Big\|_{L^2}&\lesssim e^{-ct}t^{-\frac{n}{2\delta_2}(\frac{1}{m_3}-\frac{1}{2})-\frac{s-a_3}{2\delta_2}}\|u_0\|_{H^{a_3}_{m_3}}, \label{pro3.4.3}
\end{align}
and
\small
\begin{align}
\Big\||D|^s \Big(\mathfrak{F}^{-1}\Big(e^{-t|\xi|^{2\delta_2}}\dfrac{\sin\big(t|\xi|^\sigma\big)}{|\xi|^{\sigma}}\big(1-\chi(\xi)\big)\Big) \ast u_1(x)\Big)(t,\cdot)\Big\|_{L^2}&\lesssim e^{-ct}t^{-\frac{n}{2\delta}(\frac{1}{m_3}-\frac{1}{2})-\frac{s-a_3}{2\delta_2}}\|u_1\|_{H^{[a_3-\sigma]^+}_{m_3}}. \label{pro3.4.4}
\end{align}
\normalsize
Then, combining (\ref{lemma2.5.4}), (\ref{proposition3.3.2}) and (\ref{pro3.4.2}) we may conclude (\ref{proposition3.4.2}). Combining (\ref{lemma2.5.5}), (\ref{proposition3.3.3}) and (\ref{pro3.4.3}) we may conclude (\ref{proposition3.4.3}). Combining (\ref{lemma2.5.5}), (\ref{proposition3.3.4}) and (\ref{pro3.4.4}) we may conclude (\ref{proposition3.4.4}). Summurizing, the proof of Proposition \ref{proposition3.4} is completed.
\end{proof}

\begin{proof}[Proof of Theorem \ref{theorem1.2}]
In order to prove the asymptotic profile of solutions to (\ref{equation1.1}), we may estimate
\begin{align*}
&\Big\||D|^s \Big(u(t,\cdot)- P_1\,\mathfrak{F}^{-1}\Big(e^{-\frac{1}{2}t|\xi|^{2\delta_2}}\frac{\sin\big(t|\xi|^\sigma\big)}{|\xi|^\sigma}\Big)(t,\cdot)\Big)(t,\cdot)\Big\|_{L^2} \\
&\qquad \lesssim \big\||D|^s \big(K^{\cos}_0(t,x) \ast u_0(x)\big)(t,\cdot)\big\|_{L^2}+ \big\||D|^s \big(K^{\sin}_0(t,x) \ast u_0(x)\big)(t,\cdot)\big\|_{L^2} \\
&\qquad \quad + \Big\||D|^s \Big(\Big(K_1(t,x)- \mathfrak{F}^{-1}\Big(e^{-\frac{1}{2}t|\xi|^{2\delta_2}}\frac{\sin\big(t|\xi|^\sigma\big)}{|\xi|^\sigma}\Big)\Big) \ast u_1(x)\Big)(t,\cdot)\Big\|_{L^2} \\
&\qquad \quad+ \Big\||D|^s \Big(\mathfrak{F}^{-1}\Big(e^{-\frac{1}{2}t|\xi|^{2\delta_2}}\frac{\sin\big(t|\xi|^\sigma\big)}{|\xi|^\sigma}\Big) \ast u_1(x)- P_1\,\mathfrak{F}^{-1}\Big(e^{-\frac{1}{2}t|\xi|^{2\delta_2}}\frac{\sin\big(t|\xi|^\sigma\big)}{|\xi|^\sigma}\Big)\Big)(t,\cdot)\Big\|_{L^2} \\
&\qquad =: J_1+ J_2+ J_3+ J_4.
\end{align*}
Combining (\ref{lemma2.5.1}) and (\ref{lemma2.5.4}), (\ref{lemma2.5.2}) and (\ref{lemma2.5.4}) we obtain
\begin{align*}
J_1 &\lesssim (1+t)^{-\frac{n}{4\delta_2}-\frac{s}{2\delta_2}} \|u_0\|_{\mathcal{A}^2_0}, \\ 
J_2 &\lesssim (1+t)^{-\frac{n}{4\delta_2}-\frac{s}{2\delta_2}-1+\frac{\sigma}{2\delta_2}} \|u_0\|_{\mathcal{A}^2_0},
\end{align*}
respectively. By (\ref{proposition3.4.4}) we derive
$$ J_3 \lesssim (1+t)^{-\frac{n}{2\delta_2}(\frac{1}{m}-\frac{1}{2})-1-\frac{s}{2\delta_2}+\frac{\sigma}{\delta_2}}\|u_1\|_{\mathcal{A}^2_1}. $$
To control $J_4$, we shall apply Lemma \ref{L^1.Lemma}. First, thanks to $\big|\sin\big(t|\xi|^\sigma\big)\big| \le 1$, we employ Parseval-Plancherel formula and Lemma \ref{LemmaL1normEstimate} to get
\begin{equation}
\Big\||D|^s \mathfrak{F}^{-1}\Big(e^{-\frac{1}{2}t|\xi|^{2\delta_2}}\frac{\sin\big(t|\xi|^\sigma\big)}{|\xi|^\sigma}\Big)(t,\cdot)\Big\|_{L^2}= \Big\||\xi|^{s-\sigma}e^{-\frac{1}{2}t|\xi|^{2\delta_2}}\sin\big(t|\xi|^\sigma\big)\Big\|_{L^2} \le C\,t^{-\frac{n}{4\delta_2}-\frac{s}{2\delta_2}+\frac{\sigma}{2\delta_2}}, \label{the1.2.1}
\end{equation}
where we used the condition $n> 2\sigma$ and $C$ is a suitable positive constant. Then, by Lemma \ref{L^1.Lemma} we imply immediately
$$ J_4= o\Big(t^{-\frac{n}{4\delta_2}-\frac{s}{2\delta_2}+\frac{\sigma}{2\delta_2}}\Big) \quad \text{ as } t\to \ity. $$
Therefore, from all the above estimates for $J_k$ with $k=1,2,3,4$ we may arrive at (\ref{theorem1.2.1}). Next, from (\ref{theorem1.2.1}) and (\ref{the1.2.1}) we may estimate the desired upper bound in the following way:
\begin{align*}
\big\||D|^s u(t,\cdot)\big\|_{L^2} &\le \Big\||D|^s \Big(u(t,x)- P_1\,\mathfrak{F}^{-1}\Big(e^{-\frac{1}{2}t|\xi|^{2\delta_2}}\frac{\sin\big(t|\xi|^\sigma\big)}{|\xi|^\sigma}\Big)\Big)(t,\cdot)\Big\|_{L^2} \\
&\qquad+ |P_1|\,\Big\||D|^s \mathfrak{F}^{-1}\Big(e^{-\frac{1}{2}t|\xi|^{2\delta_2}}\frac{\sin\big(t|\xi|^\sigma\big)}{|\xi|^\sigma}\Big)(t,\cdot)\Big\|_{L^2} \\
&\le C|P_1|\,t^{-\frac{n}{4\delta_2}-\frac{s}{2\delta_2}+\frac{\sigma}{2\delta_2}}+ o\Big(t^{-\frac{n}{4\delta_2}-\frac{s}{2\delta_2}+\frac{\sigma}{2\delta_2}}\Big) \\
&\le C_2\,t^{-\frac{n}{4\delta_2}-\frac{s}{2\delta_2}+\frac{\sigma}{2\delta_2}}\quad \text{ as } t\to \ity,
\end{align*}
where $C_2$ is a suitable positive constant. Moreover, to indicate the lower bound, using again Parseval-Plancherel formula and the change of variables $|\xi|= t^{-\frac{1}{2\delta_2}}|\eta|$ we have
\small
\begin{align*}
&\Big\||D|^s \mathfrak{F}^{-1}\Big(e^{-\frac{1}{2}t|\xi|^{2\delta_2}}\frac{\sin\big(t|\xi|^\sigma\big)}{|\xi|^\sigma}\Big)(t,\cdot)\Big\|_{L^2}= \Big(\int_{\R^n}|\xi|^{2(s-\sigma)}e^{-t|\xi|^{2\delta_2}}\sin^2\big(t|\xi|^\sigma\big)d\xi\Big)^{\frac{1}{2}} \\
&\qquad= \Big(\int_{\R^n}|\xi|^{2(s-\sigma)}e^{-t|\xi|^{2\delta_2}}d\xi -\frac{1}{2}\int_{\R^n}|\xi|^{2(s-\sigma)}e^{-t|\xi|^{2\delta_2}}\cos\big(2t|\xi|^\sigma\big)d\xi\Big)^{\frac{1}{2}} \\
&\qquad= t^{-\frac{n}{4\delta_2}-\frac{s}{2\delta_2}+\frac{\sigma}{2\delta_2}}\Big(\int_0^\ity |\eta|^{2(s-\sigma)+n-1}e^{-|\eta|^{2\delta_2}}d|\eta|- \frac{1}{2}\int_0^\ity |\eta|^{2(s-\sigma)+n-1}e^{-|\eta|^{2\delta_2}}\cos\big(2t^{1-\frac{\sigma}{2\delta_2}}|\eta|^\sigma\big)d|\eta|\Big)^{\frac{1}{2}}.
\end{align*}
\normalsize
Applying Lemma \ref{Rie-Les.Lemma} leads to
$$ \int_0^\ity |\eta|^{2(s-\sigma)+n-1}e^{-|\eta|^{2\delta_2}}\cos\big(2t^{1-\frac{\sigma}{2\delta_2}}|\eta|^\sigma\big)d|\eta| \to 0 \quad \text{ as } t\to \ity. $$
It follows immediately
$$ \Big\||D|^s \mathfrak{F}^{-1}\Big(e^{-\frac{1}{2}t|\xi|^{2\delta_2}}\frac{\sin\big(t|\xi|^\sigma\big)}{|\xi|^\sigma}\Big)(t,\cdot)\Big\|_{L^2} \ge C^*\,t^{-\frac{n}{4\delta_2}-\frac{s}{2\delta_2}+\frac{\sigma}{2\delta_2}} \quad \text{ as } t\to \ity, $$
where $C^*$ is a suitable positive constant. Here we notice that the integral $\int_0^\ity |\eta|^{2(s-\sigma)+n-1}e^{-|\eta|^{2\delta_2}}d|\eta|$ is a positive constant because of the condition $n> 2\sigma$. Thus, we conclude
\begin{align*}
\big\||D|^s u(t,\cdot)\big\|_{L^2}&\ge |P_1|\,\Big\||D|^s \mathfrak{F}^{-1}\Big(e^{-\frac{1}{2}t|\xi|^{2\delta_2}}\frac{\sin\big(t|\xi|^\sigma\big)}{|\xi|^\sigma}\Big)(t,\cdot)\Big\|_{L^2} \\
&\qquad- \Big\||D|^s \Big(u(t,x)- P_1\,\mathfrak{F}^{-1}\Big(e^{-\frac{1}{2}t|\xi|^{2\delta_2}}\frac{\sin\big(t|\xi|^\sigma\big)}{|\xi|^\sigma}\Big)\Big)(t,\cdot)\Big\|_{L^2} \\ 
&\ge C^*|P_1|\,t^{-\frac{n}{4\delta_2}-\frac{s}{2\delta_2}+\frac{\sigma}{2\delta_2}}- o\Big(t^{-\frac{n}{4\delta_2}-\frac{s}{2\delta_2}+\frac{\sigma}{2\delta_2}}\Big) \ge C_1\,t^{-\frac{n}{4\delta_2}-\frac{s}{2\delta_2}+\frac{\sigma}{2\delta_2}}\quad \text{ as } t\to \ity,
\end{align*}
where $C_1$ is a suitable positive constant. Summurizing, Theorem \ref{theorem1.2} is proved.
\end{proof}

\subsection{The case $a=1$ and $b=1$} \label{Sec.3.3}
In order to prove Theorem \ref{theorem1.3}, the following auxilliary results come into play.
\bmd \label{proposition3.5}
Let $s\ge 0$ and $j=0,\,1$. Let us assume $\delta_1+\delta_2>\sigma$. Then, the following estimates hold for $m \in [1,2)$:
\begin{align}
&\Big\|\partial_t^j |D|^s \Big(\Big(K^1_{0\text{\fontshape{n}\selectfont low}}(t,x)- \mathfrak{F}^{-1}\Big(e^{-t|\xi|^{2(\sigma-\delta_1)}}\chi(\xi)\Big)\Big) \ast u_0(x)\Big)(t,\cdot)\Big\|_{L^2} \nonumber \\
&\qquad \lesssim (1+t)^{-\frac{n}{2(\sigma-\delta_1)}(\frac{1}{m}-\frac{1}{2})-\frac{s}{2(\sigma-\delta_1)}-j-\frac{\sigma-2\delta_1}{\sigma-\delta_1}}\|u_0\|_{L^m}, \label{proposition3.5.1}
\end{align}
for any space dimensions $n\ge 1$, and
\begin{align}
&\Big\|\partial_t^j |D|^s \Big(\Big(K^1_{1\text{\fontshape{n}\selectfont low}}(t,x)- \mathfrak{F}^{-1}\Big(\dfrac{e^{-t|\xi|^{2(\sigma-\delta_1)}}}{|\xi|^{2\delta_1}}\chi(\xi)\Big)\Big) \ast u_1(x)\Big)(t,\cdot)\Big\|_{L^2} \\
&\qquad \lesssim (1+t)^{-\frac{n}{2(\sigma-\delta_1)}(\frac{1}{m}-\frac{1}{2})-\frac{s}{2(\sigma-\delta_1)}-j-\frac{\sigma-3\delta_1}{\sigma-\delta_1}}\|u_1\|_{L^m}, \label{proposition3.5.2}
\end{align}
for any space dimensions $n> 2m_0\delta_1$.
\emd

\begin{proof}
The proof of this proposition is similar to the proof of Proposition \ref{proposition3.1}. For this reason, we only present the steps which are different. Then, we shall repeat some of the arguments as we did in the proof of Proposition \ref{proposition3.1} to conclude the desired estimates. \\
Indeed, following the proof of Proposition \ref{proposition3.1} it is sufficient to prove the following estimate:
\begin{align}
&|\xi|^s \chi(|\xi|)\Big|\widehat{K^1_0}(t,\xi)- e^{-t|\xi|^{2(\sigma-\delta_1)}}\Big| \nonumber \\ 
&\qquad \lesssim e^{-ct|\xi|^{2(\sigma-\delta_1)}}\big(t\,|\xi|^{s+2(2\sigma-3\delta_1)+2j(\sigma-\delta_1)}+ |\xi|^{s+2(\sigma-2\delta_1)+2j(\sigma-\delta_1)}\big), \label{pro3.6.1}
\end{align}
where $c$ is a suitable positive constant. Recalling the characteristic root $\lambda_1$ we re-write as follows:
$$ -\lambda_1= \frac{1}{2}\Big(\big(|\xi|^{2\delta_1}+ |\xi|^{2\delta_2}\big)- \sqrt{\big(|\xi|^{2\delta_1}+ |\xi|^{2\delta_2}\big)^2-4|\xi|^{2\sigma}}\Big)= p\Big(1-\sqrt{1- q^2}\Big), $$
where we introduce $p:= \frac{1}{2}\big(|\xi|^{2\delta_1}+ |\xi|^{2\delta_2}\big)$ and $q:= \frac{2|\xi|^\sigma}{|\xi|^{2\delta_1}+ |\xi|^{2\delta_2}}$. It is clear that $q< 1$ for small frequencies. For this reason, applying Newton's binomial theorem gives
\begin{align*}
-\lambda_1 &= p\Big(1- \Big(1- \frac{1}{2}q^2- \frac{1}{8}q^4- o\big(q^4\big)\Big)\Big)= p\Big(\frac{1}{2}q^2+ \frac{1}{8}q^4+ o\big(q^4\big)\Big) \\
&= \frac{|\xi|^{2\sigma}}{|\xi|^{2\delta_1}+ |\xi|^{2\delta_2}}+ \frac{|\xi|^{4\sigma}}{(|\xi|^{2\delta_1}+ |\xi|^{2\delta_2})^3}+ o\Big(\frac{|\xi|^{4\sigma}}{(|\xi|^{2\delta_1}+ |\xi|^{2\delta_2})^3}\Big).
\end{align*}
Hence, a standard calculation leads to
\begin{align*}
-\lambda_1- |\xi|^{2(\sigma-\delta_1)}&= -\frac{|\xi|^{2\sigma+2(\delta_2-\delta_1)}}{|\xi|^{2\delta_1}+ |\xi|^{2\delta_2}}+ \frac{|\xi|^{4\sigma}}{(|\xi|^{2\delta_1}+ |\xi|^{2\delta_2})^3}+ o\Big(\frac{|\xi|^{4\sigma}}{(|\xi|^{2\delta_1}+ |\xi|^{2\delta_2})^3}\Big) \\ 
&= \frac{|\xi|^{2\sigma}}{|\xi|^{2\delta_1}+ |\xi|^{2\delta_2}}\Big(\frac{|\xi|^{2\sigma}}{(|\xi|^{2\delta_1}+ |\xi|^{2\delta_2})^2}- |\xi|^{2(\delta_2-\delta_1)}\Big)+ o\Big(\frac{|\xi|^{4\sigma}}{(|\xi|^{2\delta_1}+ |\xi|^{2\delta_2})^3}\Big) \\
&= \frac{|\xi|^{2\sigma}}{|\xi|^{2\delta_1}+ |\xi|^{2\delta_2}}\frac{|\xi|^{2\sigma}- |\xi|^{2(\delta_1+\delta_2)}\big(1+|\xi|^{2(\delta_2-\delta_1)}\big)^2}{(|\xi|^{2\delta_1}+ |\xi|^{2\delta_2})^2}+ o\Big(\frac{|\xi|^{4\sigma}}{(|\xi|^{2\delta_1}+ |\xi|^{2\delta_2})^3}\Big)> 0
\end{align*}
for small frequencies, where the condition $\delta_1+\delta_2>\sigma$ plays an important role. This implies immediately the following relations:
$$ \min\big\{-\lambda_1,\,|\xi|^{2(\sigma-\delta_1)}\big\}= |\xi|^{2(\sigma-\delta_1)}\quad \text{ and }\quad \Big|-\lambda_1- |\xi|^{2(\sigma-\delta_1)}\Big|\le \frac{|\xi|^{4\sigma}}{(|\xi|^{2\delta_1}+ |\xi|^{2\delta_2})^3}\le |\xi|^{2(2\sigma-3\delta_1)}. $$
By analogous arguments as in the proof of Proposition \ref{proposition3.1} we may arrive at (\ref{pro3.6.1}). Therefore, the proof of Proposition \ref{proposition3.5} is completed.
\end{proof}

\noindent Following the proof of Proposition \ref{proposition3.2} we may conclude the following result by using Proposition \ref{proposition3.5} and Lemma \ref{lemma2.6}.
\bmd \label{proposition3.6}
Let $s\ge 0$ and $j=0,\,1$. Let us assume $\delta_1+\delta_2>\sigma$. Then, the following estimates hold for $m \in [1,2)$:
\begin{align}
&\Big\|\partial_t^j |D|^s \Big(\Big(K^1_0(t,x)- \mathfrak{F}^{-1}\Big(e^{-t|\xi|^{2(\sigma-\delta_1)}}\Big)\Big) \ast u_0(x)\Big)(t,\cdot)\Big\|_{L^2} \nonumber \\
&\qquad \lesssim (1+t)^{-\frac{n}{2(\sigma-\delta_1)}(\frac{1}{m}-\frac{1}{2})-\frac{s}{2(\sigma-\delta_1)}-j-\frac{\sigma-2\delta_1}{\sigma-\delta_1}}\|u_0\|_{L^m}+ e^{-ct}t^{-\frac{n}{2(\sigma-\delta_2)}(\frac{1}{m_1}-\frac{1}{2})-\frac{s-a_1}{2(\sigma-\delta_2)}-j}\|u_0\|_{H^{a_1}_{m_1}} \nonumber \\
&\qquad \quad +e^{-ct}t^{-\frac{n}{2\delta_2}(\frac{1}{m_2}-\frac{1}{2})-\frac{s-a_2}{2\delta_2}-j+\frac{2\delta_2-\sigma}{\delta_2}}\|u_0\|_{H^{a_2}_{m_2}} +e^{-ct}t^{-\frac{n}{2(\sigma-\delta_1)}(\frac{1}{m_3}-\frac{1}{2})-\frac{s-a_3}{2(\sigma-\delta_1)}-j}\|u_0\|_{H^{a_3}_{m_3}} \label{proposition3.6.1}
\end{align}
for any space dimensions $n\ge 1$, and
\small
\begin{align}
&\Big\|\partial_t^j |D|^s \Big(\Big(K^1_1(t,x)- \mathfrak{F}^{-1}\Big(\dfrac{e^{-t|\xi|^{2(\sigma-\delta_1)}}}{|\xi|^{2\delta_1}}\Big)\Big) \ast u_1(x)\Big)(t,\cdot)\Big\|_{L^2} \nonumber \\
&\qquad \lesssim (1+t)^{-\frac{n}{2(\sigma-\delta_1)}(\frac{1}{m}-\frac{1}{2})-\frac{s}{2(\sigma-\delta_1)}-j-\frac{\sigma-3\delta_1}{\sigma-\delta_1}}\|u_1\|_{L^m}+ e^{-ct}t^{-\frac{n}{2(\sigma-\delta_2)}(\frac{1}{m_1}-\frac{1}{2})-\frac{s-a_1}{2(\sigma-\delta_2)}-j}\|u_1\|_{H^{[a_1-2\delta_2]^+}_{m_1}} \nonumber \\
&\qquad \quad + e^{-ct}t^{-\frac{n}{2\delta_2}(\frac{1}{m_2}-\frac{1}{2})-\frac{s-a_2}{2\delta_2}-j}\|u_1\|_{H^{[a_2-2\delta_2]^+}_{m_2}} +e^{-ct}t^{-\frac{n}{2(\sigma-\delta_1)}(\frac{1}{m_3}-\frac{1}{2})-\frac{s-a_3}{2(\sigma-\delta_1)}-j}\|u_1\|_{H^{[a_3-2\delta_1]^+}_{m_3}} \label{proposition3.6.2}
\end{align}
\normalsize
for any space dimensions $n> 2m_0\delta_1$. Here $a_1,\,a_2,\,a_3 \ge 0$, $m_1,\,m_2,\,m_3 \in [1,2]$ and $c$ is a suitable positive constant.
\emd

\begin{proof}[Proof of Theorem \ref{theorem1.3}]
In order to prove the asymptotic profile of solutions to (\ref{equation1.1}), we will follow the proof of Theorem \ref{theorem1.1}. At first, as in the proof of Theorem \ref{theorem1.1} we recall the following estimates:
\begin{align*}
&\Big\|\partial_t^j |D|^s \Big(u(t,\cdot)- P_1\,\mathfrak{F}^{-1}\Big(\dfrac{e^{-t|\xi|^{2(\sigma-\delta_1)}}}{|\xi|^{2\delta_1}}\Big)(t,\cdot)\Big)(t,\cdot)\Big\|_{L^2} \\
&\quad \lesssim \big\|\partial_t^j |D|^s \big(K^1_0(t,x) \ast u_0(x)\big)(t,\cdot)\big\|_{L^2}+ \big\|\partial_t^j |D|^s \big(K^2_0(t,x) \ast u_0(x)\big)(t,\cdot)\big\|_{L^2} \\
&\qquad+ \big\||D|^s \big(K^2_1(t,x) \ast u_1(x)\big)(t,\cdot)\big\|_{L^2}+ \Big\|\partial_t^j|D|^s \Big(\Big(K^1_1(t,x)- \mathfrak{F}^{-1}\Big(\frac{e^{-t|\xi|^{2(\sigma-\delta_1)}}}{|\xi|^{2\delta_1}}\Big)\Big) \ast u_1(x)\Big)(t,\cdot)\Big\|_{L^2} \\
&\qquad+ \Big\|\partial_t^j |D|^s \Big(\mathfrak{F}^{-1}\Big(\frac{e^{-t|\xi|^{2(\sigma-\delta_1)}}}{|\xi|^{2\delta_1}}\Big) \ast u_1(x)- P_1\,\mathfrak{F}^{-1}\Big(\dfrac{e^{-t|\xi|^{2(\sigma-\delta_1)}}}{|\xi|^{2\delta_1}}\Big)\Big)(t,\cdot)\Big\|_{L^2} \\
&\quad =: I_1+ I_2+ I_3+ I_4+ I_5.
\end{align*}
Combining (\ref{lemma2.6.1}) and (\ref{lemma2.6.5}), (\ref{lemma2.6.2}) and (\ref{lemma2.6.5}), (\ref{lemma2.6.4}) and (\ref{lemma2.6.6}) we derive
\begin{align*}
I_1 &\lesssim (1+t)^{-\frac{n}{4(\sigma-\delta_1)}- \frac{s}{2(\sigma-\delta_1)}-j} \|u_0\|_{\mathcal{A}^3_0}, \\ 
I_2 &\lesssim (1+t)^{-\frac{n}{4\delta_1}-\frac{s}{2\delta_1}-j- \frac{\sigma-2\delta_1}{\delta_1}} \|u_0\|_{\mathcal{A}^3_0}, \\
I_3 &\lesssim (1+t)^{-\frac{n}{4\delta_1}- \frac{s}{2\delta_1}-j+1} \|u_1\|_{\mathcal{A}^3_1},
\end{align*}
respectively. By (\ref{proposition3.6.2}) we arrive at
$$ I_4 \lesssim (1+t)^{-\frac{n}{2(\sigma-\delta_1)}(\frac{1}{m}-\frac{1}{2})-\frac{s}{2(\sigma-\delta_1)}-j-\frac{\sigma-3\delta_1}{\sigma-\delta_1}}\|u_1\|_{\mathcal{A}^3_1}. $$
Then, we shall repeat some of the arguments as we did in the proof of Theorem \ref{theorem1.1} to conclude the desired estimates.
Summurizing, Theorem \ref{theorem1.3} is proved.
\end{proof}


\noindent \textbf{Acknowledgments}\medskip

\noindent The PhD study of the second author is supported by Vietnamese Government's Scholarship (Grant number: 2015/911).\medskip

\noindent\textbf{Appendix A}\medskip

\noindent \textit{A.1. Useful lemmas}
\bbd \label{LemmaL1normEstimate}
Let $n\ge1$, $c>0$, $\alpha >0$ and $\beta \in \R$ satisfy $n+\beta >0$. 
The following estimates hold for $t>0${\rm :}
$$ \int_{|\xi|\le1} |\xi|^\beta e^{-c|\xi|^\alpha t}d\xi \lesssim (1+t)^{-\frac{n+\beta}{\alpha}} \quad \text{ and } \quad \int_{|\xi|\ge1} |\xi|^\beta e^{-c|\xi|^\alpha t}d\xi \lesssim t^{-\frac{n+\beta}{\alpha}}. $$
\ebd
The proof of this lemma can be found in \cite{DaoMichihisa}.

\bbd \label{L^1.Lemma}
Let $a\ge 0$. Let us assume $v= v(x) \in L^1$ and $\phi=\phi(t,x)$ be a smooth function satisfying
$$ \big\||D|^a \phi(t,\cdot)\big\|_{L^2} \lesssim t^{-\alpha} \quad \text{ and }\quad  \big\||D|^{a+1} \phi(t,\cdot)\big\|_{L^2} \lesssim t^{-\alpha-\beta}, $$
for some positive constants $\alpha,\,\beta>0$. Then, the following estimate holds:
$$ \Big\||D|^a \Big(\phi(t,x) \ast v(x)- \Big(\int_{\R^n}v(y)\,dy\Big)\phi(t,x)\Big)(t,\cdot) \Big\|_{L^2}= o\big(t^{-\alpha}\big) \quad \text{ as }t \to \ity, $$
for all space dimensions $n\ge 1$.
\ebd
One can be found the proof of this lemma in \cite{IkehataTakeda}.

\bbd[A variant version of Riemann-Lebesgue theorem] \label{Rie-Les.Lemma}
If $f= f(r) \in L^1$ and \text{\fontshape{n}\selectfont supp}$f \in (0,\ity)$, then it holds:
$$ \int_0^\ity f(r)e^{-zr}dr \to 0 \quad \text{ as }|z|\to \ity \text{ within the half-plane }\text{\fontshape{n}\selectfont Re}z \ge 0, $$
that is,
$$ \int_0^\ity f(r)\cos(r\tau)dr \to 0 \quad \text{ and }\quad \int_0^\ity f(r)\sin(r\tau)dr \to 0 \quad \text{ as } \tau \to \ity. $$
\ebd


\end{document}